\documentclass{amsart}
\usepackage[T1]{fontenc}
\usepackage{tgtermes}
\usepackage{amsmath, amsthm, amsfonts, amssymb}
\usepackage{mathrsfs}
\usepackage{graphicx}
\usepackage{hyperref}
\usepackage{url}
\usepackage[all]{xy}
\usepackage{enumerate}
\usepackage{fancyhdr}
\usepackage{xspace}
\usepackage{verbatim}

\usepackage{tikz,xcolor}
\usetikzlibrary{arrows,decorations.pathmorphing,backgrounds,fit,automata}
\usetikzlibrary{shapes.geometric,shapes.arrows,calc}
\usetikzlibrary{positioning}
\theoremstyle{definition}
\newtheorem*{ack}{Acknowledgements}
\newtheorem{defn}{Definition}[section]
\newtheorem{cond}[defn]{Condition}
\newtheorem{eg}[defn]{Example}

\theoremstyle{remark}

\theoremstyle{plain}

\newtheorem{cor}[defn]{Corollary}

\newtheorem{lem}[defn]{Lemma}

\newtheorem{prop}[defn]{Proposition}

\newtheorem{thm}[defn]{Theorem}

\numberwithin{equation}{section}
\numberwithin{figure}{section}

\newcommand{\abs}[1]{\left\vert#1\right\vert}
\newcommand{\set}[1]{\left\{#1\right\}}
\newcommand{\gen}[1]{\left\langle#1\right\rangle}

\newcommand{\round}[1]{{\ooalign{\hfil\raise .10ex\hbox{\scriptsize#1}\hfil\crcr\mathhexbox20D}}}
\newcommand{\recur}[1]{\overline{#1}}

\newcommand{\B}[1]{\ensuremath{\mathbb{#1}}}

\newcommand{\C}[1]{\ensuremath{\mathcal{#1}}}

\newcommand{\alb}{\mathsf{X}}
\newcommand{\automaton}{\mathsf{A}}
\newcommand{\crits}{\mathcal{C}}
\newcommand{\pcrits}{\mathcal{P}}
\newcommand{\nucl}{\mathcal{N}}
\newcommand{\lims}{\mathcal{J}_G}
\newcommand{\sss}{\mathcal{L}}
\newcommand{\tile}{\mathcal{T}}
\newcommand{\shift}{\mathsf{s}}

\def\fr{{fi\-ni\-te\-ly ram\-i\-fied}}

\def\pcf{{\mbox{p.c.f.}\xspace}}
\def\pcfp{{\mbox{p.c.f}\@. }}
\def\s-s{self-similar}

\DeclareMathOperator{\Aut}{Aut}

\allowdisplaybreaks[4]
\title{From self-similar structures to self-similar groups}
\author[D.~J.~Kelleher]{Daniel J.~Kelleher}
\address[D.~J.~Kelleher]{Department of Mathematics, University of Connecticut, Storrs, CT 06269, USA}
\email{kelleher@math.uconn.edu}
\urladdr{\url{http://www.math.uconn.edu/~kelleher/}}

\author[B.~A.~Steinhurst]{Benjamin A.~Steinhurst}
\address[B.~A.~Steinhurst]{Department of Mathematics, Cornell University, Ithaca, NY 14853, USA}
\email{steinhurst@math.uconn.edu}
\urladdr{\url{http://www.math.cornell.edu/~steinhurst/}}

\author[C.-M.~M.~Wong]{Chuen-Ming~M.~Wong}
\address[C.-M.~M.~Wong]{Department of Mathematics, Princeton University, Princeton, NJ 08544, USA}
\email{chuenw@princeton.edu}
\urladdr{\url{http://www.princeton.edu/~chuenw/}}

\thanks{Authors supported in part by the National Science Foundation through grant DMS-0505622.}
\subjclass[2000]{Primary 20E08; Secondary 20F65}
\keywords{Fractals, \s-sity, limit space, \pcf,}
\begin{document}
%
%
\begin{abstract}
We explore the relationship between limit spaces of contracting self-similar groups and self-similar structures. We give the condition on a contracting group such that its limit space admits a self-similar structure, and also the condition such that this self-similar structure is p.c.f. We then give the necessary and sufficient condition on a p.c.f. self-similar structure such that there exists a contracting group whose limit space has an isomorphic self-similar structure; in this case, we provide a construction that produces such a contracting group. Finally, we illustrate our results with several examples.
\end{abstract}

\maketitle

\section{Introduction}

The theory of \s-s groups developed as a part of geometric group theory in the last decades. In this theory, many exotic groups (such as groups of intermediate growth, non-elementary amenable groups, and infinite finitely generated torsion groups) could be easily described by their actions on a rooted tree \cite{Gri80, Gri84, GS83, BGS03}. More recently, a close relationship between the theory of \s-s groups and fractals has been discovered and studied \cite{BGN03, Nek05, NT08}. This survey aims to clarify this relationship by closely examining the correspondence between \s-s groups and \s-s structures on fractals.

The standard reference on the theory of \s-s groups is \cite{Nek05}. A \emph{\s-s group} is an automorphism group acting on the rooted tree in a recursive manner. Every \emph{contracting} \s-s group $G$ (see Definition~\ref{defn:contracting}) induces an \emph{asymptotic equivalence relation} on the \emph{boundary} of the rooted tree (i.e.\ the space of left-finite words), and the quotient space of the boundary of the rooted tree by this equivalence relation is called the \emph{limit space} $\lims$ of the \s-s group. $G$ does not act on $\lims$, but rather contains information about the adjacency of cells of $\lims$ and describes its fractal-like properties.

Of these properties, we are most interested in the \s-sity of the limit space; to examine this, we employ the notion of \emph{\s-s structures}, which is fundamental in analysis on fractals. Self-similar structures have been extensively studied in \cite{Kig01}. A \s-s structure is a finite set of injections from a compact space $K$ to itself, such that $K$ is covered by the union of the images. By repeatedly applying these injections, each point in $K$ can then be given some \emph{addresses} in the \emph{code space}, which can be identified with the boundary of a rooted tree; in this way, $K$ can also be viewed as a quotient space of the boundary of the rooted tree. Many well-known fractals can be given a \s-s structure.

To capture the \s-sity of the limit space of a contracting \s-s group, we only consider the \s-s structure on the limit space where the two quotient maps mentioned above are the same. Such a \s-s structure may not exist (see Example~\ref{eg:adding_machine}), but if it does, then the action of the group can aid the development of analytic properties of the limit space (see \cite{NT08}). A natural question arises: Under what conditions on the \s-s group does this \s-s structure exist? We attempt to give an answer in Theorem~\ref{thm:quasi-monocarpic}.

In the second half of the paper, we focus on \emph{post-critically finite} (\pcf) \s-s structures on limit spaces. A \pcf\ fractal is one where the cells of the fractal only intersect at finitely many points, and these points have finitely many addresses. This is a broad class of fractals on which the methods to develop a Laplacian, as well as an analogous Gaussian process, are known. For more details on \pcf\ \s-s structures and their analytic properties, see \cite{Kig93, Kig01, Str06}. We investigate the conditions for which the \s-s structure on a limit space is \pcf\ (see Theorem~\ref{thm:pcf_pcf}). We also attempt to answer the inverse question: Given a \pcf\ \s-s structure, can we find a contracting \s-s group whose limit space has an isomorphic \s-s structure? In Section~\ref{sec:construction}, we shall identify the condition on the \s-s structure (Theorem~\ref{thm:construct_iff_s_exists}) such that the answer is affirmative, and attempt to directly construct a contracting \s-s group in this case.

In practice, we see that certain fractals, equipped with any \s-s structure, cannot arise as the limit space of a contracting action. This includes both non-\pcf\ fractals (such as the diamond fractal) and \pcf\ fractals (such as the Linstr{\o}m snowflake). For some fractals (including the Sierpi\'nski gasket and the pentakun), a contracting group can be found only for some (and not all) \s-s structures on the space. Also, we shall exhibit that two contracting groups that are not isomorphic can have the same limit space with the same \s-s structure in Example~\ref{eg:alt_grigorchuk}.

This paper is organized as follows. We begin with a brief review of the basic definitions of \s-s groups and limit spaces in Section~\ref{sec:preliminaries}. Section~\ref{sec:sss_on_limit_spaces} makes precise the notion of the \s-s structure on a limit space, and gives the condition on the contracting group that ensures the existence of this induced \s-s structure. Building on the work in \cite{BN03, NT08}, we discuss \pcf\ \s-s structures on limit spaces in Section~\ref{sec:limit_space_pcf_action}, and in particular give the condition for the \s-s structure on a limit space to be \pcfp In Section~\ref{sec:construction}, we address the inverse problem: We detail a construction that, given a \pcf \ \s-s structure with a certain necessary condition, produces a contracting \s-s group whose limit space has an isomorphic \s-s structure. Finally, Section~\ref{sec:examples} is a compilation of examples illustrating the findings of this paper.

\begin{ack}
The authors wish to thank Alexander Teplyaev for his guidance on the direction of research, and Robert Strichartz and Volodymyr Nekrashevych for their advice and observations, which were very helpful in writing this paper.
\end{ack}

\section{Preliminaries}
\label{sec:preliminaries}

We begin by reviewing the basic definitions in the theory of self-similar groups. For more details, see \cite{BGN03, Nek05, NT08}.

Let $\alb$ be a finite set, called the \emph{alphabet}. We denote by $\alb^n$ all \emph{finite words} $w = x_n \ldots x_2 x_1$ of length $n$ over $\alb$, where $x_i \in \alb$. The length of the word $w$ is denoted by $\abs{w}$. The set of all finite words, including the empty word $\varnothing$, is denoted by $\alb^* = \bigcup_{n=0}^\infty \alb^n$. The set $\alb^*$ has a natural structure of a rooted tree with the root $\varnothing$, where a word $w \in \alb^*$ is connected by an edge to each of the words of the form $xw$, where $x \in \alb$.

Consider the set $\alb^{-\omega} = \set{\ldots x_2 x_1: x_i \in \alb}$ of all \emph{left-infinite words}. There is a natural topology on the disjoint union $\alb^* \sqcup \alb^{-\omega}$ given by the basis consisting of open sets of the form $\alb^* w \sqcup \alb^{-\omega} w$, where $\alb^* w$ and $\alb^{-\omega} w$ denote the sets of words ending by the finite word $w$. In this topology, $\alb^{-\omega}$ is homeomorphic to the countable product of the discrete set $\alb$, and therefore to the Cantor set.

For all $w \in \alb^*$, we identify $w$ with the map $\tilde{w}: \alb^{-\omega} \to \alb^{-\omega}$, defined by $\tilde{w} (\ldots x_2 x_1) = \ldots x_2 x_1 w$. We also define the shift map $\sigma: \alb^{-\omega} \to \alb^{-\omega}$ by $\sigma (\ldots x_2 x_1) = \ldots x_3 x_2$.

An \emph{automorphism} of the rooted tree $\alb^*$ is a permutation of $\alb^*$ that fixes $\varnothing$ and preserves adjacency of the vertices. The group of all automorphisms of $\alb^*$ is denoted by $\Aut \alb^*$. We shall denote the identity automorphism by $1$. Every automorphism $g \in \Aut \alb^*$ preserves the \emph{levels}, so that $\abs{g(w)} = \abs{w}$ for every $w \in \alb^*$.

Let $g \in \Aut \alb^*$. If we identify the first level $\alb^1$ of the rooted tree with $\alb$, then the restriction of $g$ to $\alb^1$ is a permutation of $\alb$, which will be called the \emph{root permutation} of $g$ and denoted $\sigma_g$. For every $x \in \alb$, if we identify both the subtrees $g\alb^*$ and $\sigma_g (x) \alb^*$, then the restriction of $g$ to $x \alb^*$ is another automorphism of $\alb^*$, called the \emph{restriction} of $g$ at $x$ and denoted $g|_x$. Then we can write $g (xw) = \sigma_g (x) g|_x (w)$ for all $x \in \alb$ and $w \in \alb^*$.

More generally, for each $v \in \alb^*$ we identify the subtrees $v \alb^*$ and $g(v) \alb^*$ and write $g (vw) = g (v) g|_v (w)$. (We define $g|_\varnothing = g$.) Then we have the basic identities
\[
g|_{v_1v_2} = g|_{v_1}|_{v_2}\text{,}
\]
\[
(g_1 g_2) |_v = g_1 |_{g_2 (v)} g_2 |_v\text{.}
\]
This gives us
\[
g_1 g_2 (xw) = g_1 g_2 (x) g_1|_{g_2(x)} g_2|_x (w)\text{.}
\]
We shall also use the ``wreath recursion'' notation to express this. If $\alb = \set{0, \ldots, k-1}$, then
\[
g = \sigma_g (g|_0, g|_1, \ldots, g|_{k-1})\text{,}
\]
\[ g_1 g_2 = \sigma_{g_1} \sigma_{g_2} (g_1|_{\sigma_{g_2} (0)} g_2|_0, \ldots, g_1|_{\sigma_{g_2} (k-1)} g_2|_{k-1})\text{.}
\]

\begin{defn}
A faithful action of a group $G$ on the rooted tree $\alb^*$ is said to be \emph{\s-s}, or \emph{state-closed}, if for every $g \in G$ and every $x \in \alb$ there exist $h \in G$ and $y \in \alb$ such that $g (xw) = y h(w)$ for every $w \in \alb^*$.
\end{defn}

We will denote such an action as $(G, \alb)$. If $g (xw) = y h(w)$, then obviously $y = g(x)$ and $h = g|_x$. We will also write $g \cdot x = y \cdot h$. Given a faithful action of $G$ on $\alb^*$, there is a natural isomorphism between $G$ and a subgroup of $\Aut \alb^*$, with which it will be identified. Thus, we will also use the terms \emph{\s-s subgroup of $\Aut \alb^*$} and \emph{\s-s automorphism group} to describe a \s-s action.

A set $S \subset \Aut \alb^*$ of automorphisms is said to be \emph{state-closed} if the restriction of every $g \in S$ to every $x \in \alb$ is also in $S$. If every element of $S$ has finite order, then the group $G = \langle S \rangle$ is \s-s.

The notion of a \s-s action is closely related to that of an automaton \cite{ECH+92, Wol02}.

\begin{defn}
An \emph{automaton} $\automaton$ over the alphabet $\alb$ is a set of \emph{internal states}, also denoted $\automaton$, together with a map $\tau: \automaton \times \alb \to \alb \times \automaton$.
\end{defn}

An automaton is \emph{finite} if and only if its set of internal states is finite. If $\tau (q, x) = (y, p)$, we will also write formally $q \cdot x = y \cdot p$. For every state $q \in \automaton$, we can define the \emph{action of the state $q$} on all finite words $w = x_n \ldots x_2 x_1$, by processing the letters one by one: it reads the first letter $x$ of $w$, outputs the letter $p = q (x)$, goes to the state $y = q|_x$ and goes on to read the next letter. At the end it will give as output some word $q (w)$ where $\abs{q (w)} = \abs{w}$, and stop at some state of $\automaton$.

An automaton $\automaton$ is often represented by its \emph{Moore diagram}, which is a directed graph with the set of vertices $\automaton$, in which for every $x \in \alb$ and every $q \in \automaton$, there is an arrow from $q$ to $q|_x$ labeled $(x, q(x))$. Then for $q \in \automaton$ and $w = x_n \ldots x_2 x_1 \in \alb^*$, the image $q (w)$ under the action of the state $q$ can be found by finding a path in the Moore diagram which starts at $q$ with consecutive labels of the form $(x_1, y_1), (x_2, y_2), \ldots, (x_n, y_n)$; then $q(w) = y_n \ldots y_2 y_1$.

The relationship between \s-s actions and automatons is illustrated below.

\begin{defn}
Let $(G, \alb)$ be a \s-s action. An automaton $\automaton$ is said to be the \emph{complete automaton} of $G$ if its set of internal states is $G$ and the action of the states coincides with the action of $G$.
\end{defn}

It is routine to prove that there is a one-to-one correspondence between the \s-s actions and their complete automatons. Henceforth, we will identify a \s-s action with its complete automaton, and denote a state of the automaton by its corresponding group element $g$.

\begin{defn}
\label{defn:contracting}
A \s-s action is said to be \emph{contracting} if there exists a finite set $\nucl \subset G$ such that for every $g \in G$, there exists $k \in \B{N}$ such that $g|_w \in \nucl$ whenever $\abs{w} \geq k$. The smallest such $\nucl$ is called the \emph{nucleus} of $G$.
\end{defn}

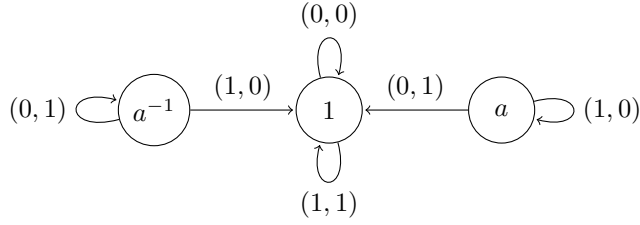
\begin{figure}
\begin{tikzpicture}[->, shorten >=1pt, node distance=1.4cm]
\node[state] at (0,0)      (1)   {$1$};
\node[state] [left=of 1]   (a-1) {$a^{-1}$};
\node[state] [right=of 1]  (a)   {$a$};
\path
(a-1) edge              node [above] {$(1,0)$} (1)
     edge [loop left=28]  node         {$(0,1)$} (a-1)
(a) edge              node [above] {$(0,1)$} (1)
     edge [loop right=28] node         {$(1,0)$} (a)
(1)  edge [loop above=28] node [above] {$(0,0)$} (1)
     edge [loop below=28] node [below] {$(1,1)$} (1);
\end{tikzpicture}
\caption{Nucleus of the binary adding machine}
\label{fig:adding_machine}
\end{figure}

Figure~\ref{fig:adding_machine} shows the Moore diagram of the nucleus of the binary adding machine, one of the simplest contracting \s-s groups. For more on this group, see Example~\ref{eg:adding_machine}.

Contracting self-similar actions have an associated topological space, which we describe below.

\begin{defn}
Let $(G, \alb)$ be a contracting action. Two left-infinite words $\ldots x_2 x_1, \ldots y_2 y_1 \in \alb^{-\omega}$ are said to be \emph{asymptotically equivalent} if there exists a sequence $\set{g_k}_{k \geq 1}$ of group elements, taking only finitely many values, such that $g (x_k \ldots x_2 x_1) = y_k \ldots y_2 y_1$ for every $k \geq 1$. The quotient space of $\alb^{-\omega}$ by the asymptotic equivalence relation is called the \emph{limit space} of the action and is denoted $\lims$; we denote the quotient map by $p: \alb^{-\omega} \to \lims$.
\end{defn}

We shall use the following more useful characterization of asymptotic equivalence.

\begin{thm}[\cite{Nek05} Theorem 3.6.3]
Let $(G, \alb)$ be a contracting action. Two left-infinite words $\ldots x_2 x_1, \ldots y_2 y_1 \in \alb^{-\omega}$ are asymptotically equivalent if and only if there exists a left-infinite path $\ldots e_2 e_1$ in the Moore diagram of the nucleus such that the edge $e_n$ is labeled by $(x_n, y_n)$.

The topological space $\lims$ is compact, metrizable and has topological dimension less than the size of the nucleus.
\end{thm}

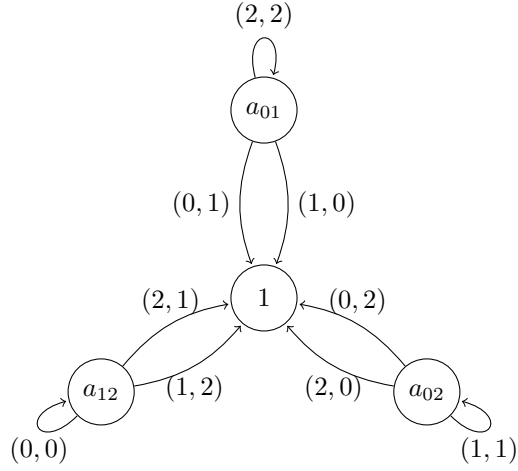
\begin{figure}
\begin{tikzpicture}[->, shorten >=1pt, node distance=2cm]
\node[state] at (0,0)     (1)  {$1$};
\node[state] at (90:2.5)  (01) {$a_{01}$};
\node[state] at (210:2.5) (12) {$a_{12}$};
\node[state] at (330:2.5) (02) {$a_{02}$};
\path
(01) edge [loop above]          node [above] {$(2,2)$} ()
     edge [bend right=20]       node [left]  {$(0,1)$} (1)
     edge [bend left=20]        node [right] {$(1,0)$} (1)
(12) edge [in=195,out=225,loop] node [below] {$(0,0)$} ()
     edge [bend right=20]       node [below] {$(1,2)$} (1)
     edge [bend left=20]        node [above] {$(2,1)$} (1)
(02) edge [in=345,out=315,loop] node [below] {$(1,1)$} ()
     edge [bend right=20]       node [above] {$(0,2)$} (1)
     edge [bend left=20]        node [below] {$(2,0)$} (1);
\end{tikzpicture}
\caption{Nucleus of the 3-peg Hanoi Towers Group}
\label{fig:hanoi_towers}
\end{figure}

Figure~\ref{fig:hanoi_towers} shows the Moore diagram of the nucleus of the 3-peg Hanoi Towers Group, which is a well-known contracting \s-s group. Its limit space has been shown in \cite{GS06, GS08} to be the Sierpi\'nski gasket. For more about this group, see Example~\ref{eg:sierp_gasket}

We now mention a property of the quotient map $p$, and the definition of the induced shift map $\shift$.

\begin{prop}
\label{prop:lim_sp_eq_rel}
Let $\lims$ be the limit space of a contracting action $(G, \alb)$ with the quotient map $p: \alb^{-\omega} \to \lims$. Then $p (\ldots x_2 x_1) = p (\ldots y_2 y_1)$ implies that $p (\ldots x_{n+1} x_n) = p (\ldots y_{n+1} y_n)$ for all $n \in \B{Z}^+$. In particular, the induced shift map $\shift: \lims \to \lims$ defined by $\shift \circ p = p \circ \sigma$ is well-defined.
\end{prop}

\begin{proof}
Since $p(\ldots x_2 x_1) = p(\ldots y_2 y_1)$, the left-infinite words $\ldots x_2 x_1$ and $\ldots y_2 y_1$ are asymptotically equivalent. Thus, there exists a left-infinite path $\ldots e_2 e_1$ in the nucleus where the label of the edge $e_k$ is $(x_k, y_k)$. Then the left-infinite path $\ldots e_{n+1} e_n$ gives the asymptotic equivalence between $\ldots x_{n+1} x_n$ and $\ldots y_{n+1} y_n$, and so $p (\ldots x_{n+1} x_n) = p (\ldots y_{n+1} y_n)$.
\end{proof}

Finally, we introduce the notion of a tile.

\begin{defn}
Let $\lims$ be the limit space of a contracting action. For each $w \in \alb^*$ such that $\abs{w} = n$, the \emph{$n^{\text{th}}$ level tile} $\tile_{w}$ is defined as the image $p (\alb^{-\omega} w)$ in $\lims$.
\end{defn}

\section{Self-similar structures on limit spaces}
\label{sec:sss_on_limit_spaces}

An important aspect of a fractal is its \s-sity. To make the adjective ``\s-s'' precise, we adopt the following definition:

\begin{defn}
Let $K$ be a compact metrizable topological space, and let $F_i: K \to K$ be a continuous injection for each $i \in \alb$, then the system $\sss = (K, \alb, \set{F_i}_{i \in \alb})$ is said to be a \emph{\s-s structure on $K$} if there exists a continuous surjection $\pi: \alb^{-\omega} \to K$ such that the relation $F_i \circ \pi = \pi \circ i$ holds, where $i(w) = wi$ for all $i\in \alb$. In this case, for each $w = x_n \ldots x_2 x_1 \in \alb^*$, we define $F_w$ by $F_w = F_{x_1} F_{x_2} \circ \ldots \circ F_{w_n}$, and the \emph{$n^{\text{th}}$ level cell} $K_w$ to be $K_w = F_w (K)$.
\end{defn}

\begin{figure}
\begin{center}
\includegraphics[width=0.4\linewidth]{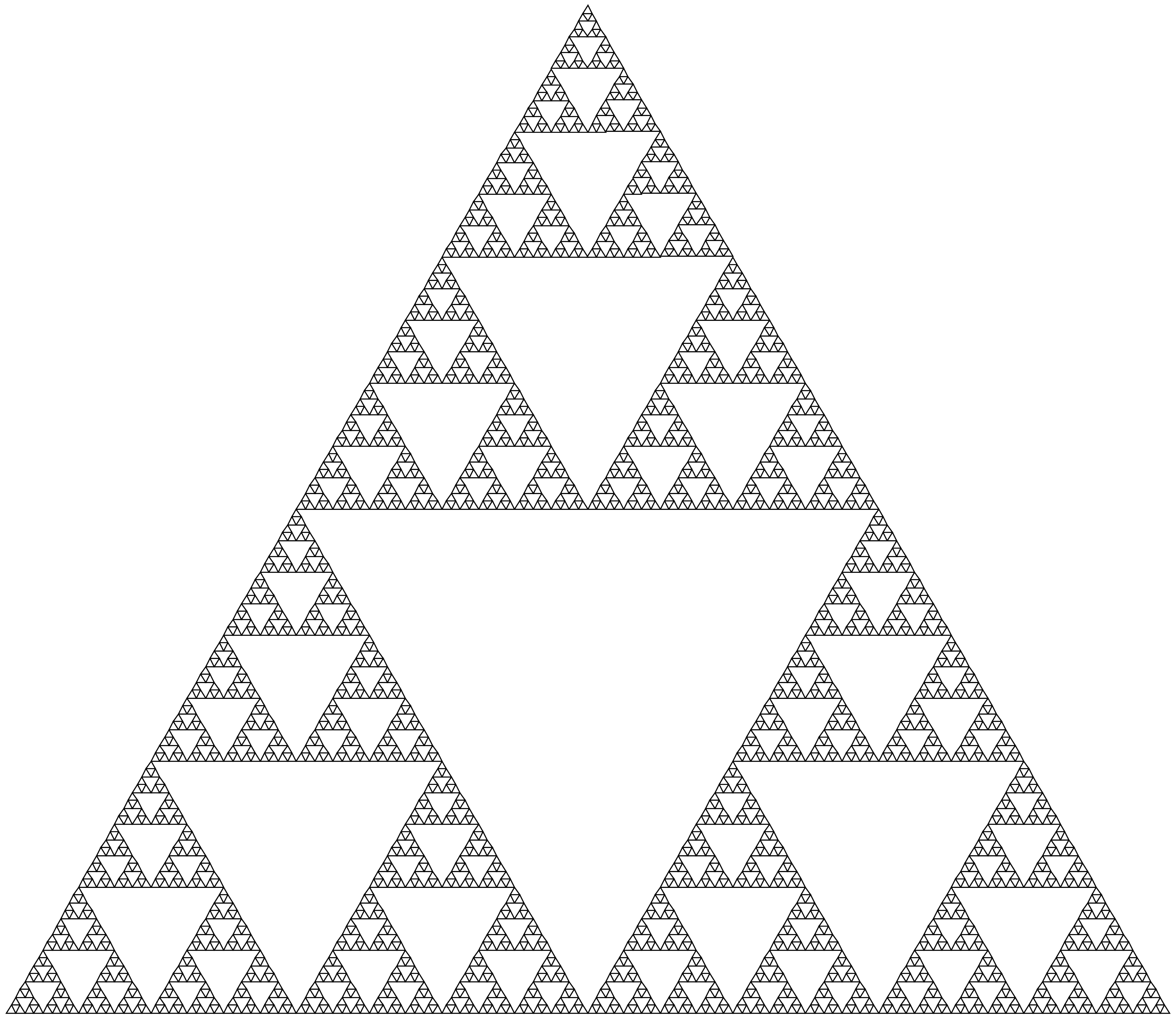}\hspace{12pt}
\includegraphics[width=0.4\linewidth]{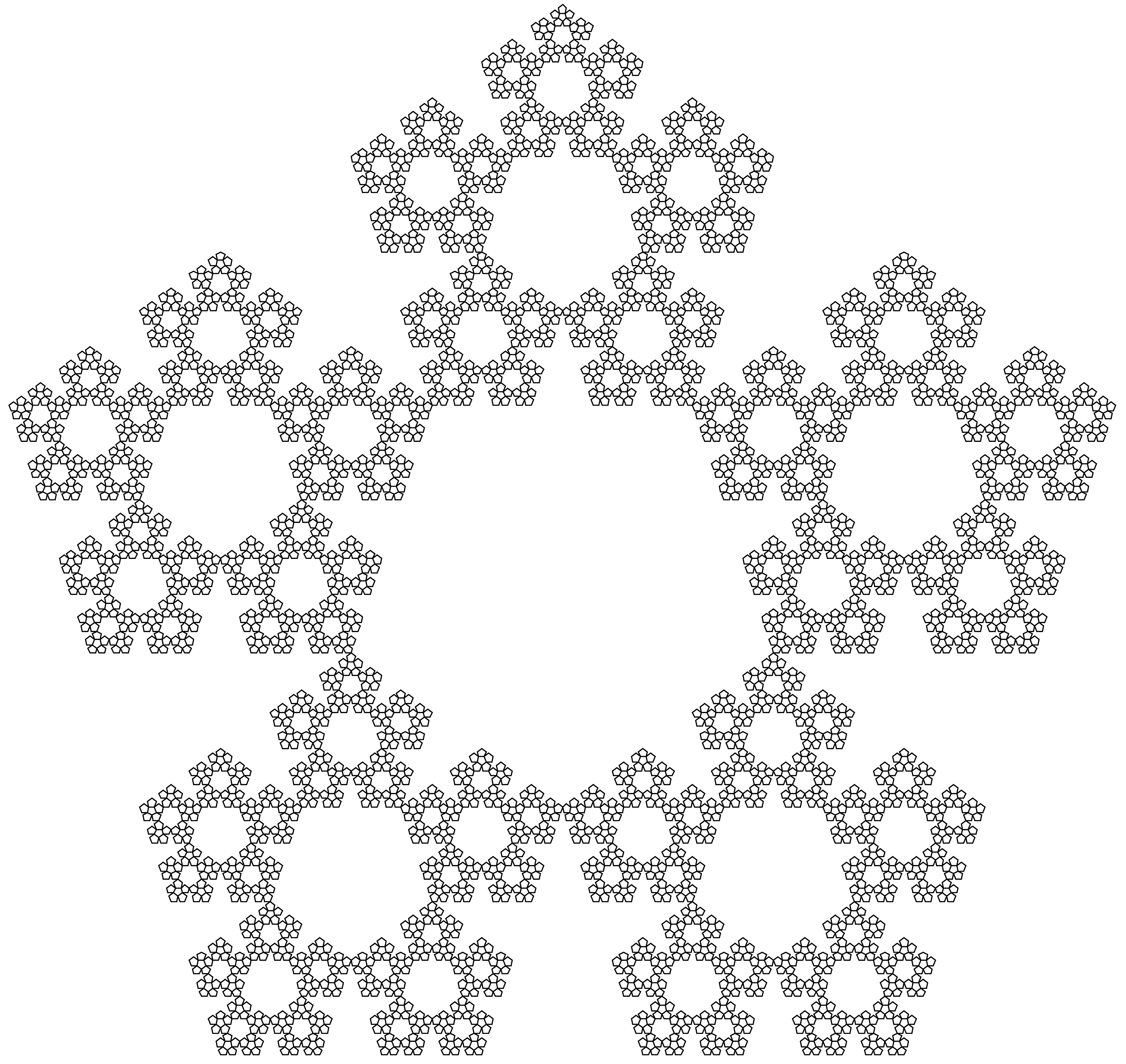}
\end{center}
\caption{Sierp\'nski gasket and pentakun}
\label{fig:sierp_gasket_pentakun}
\end{figure}

Many well-known fractals, such as the Sierpi\'nski gasket and the pentakun (see Figure~\ref{fig:sierp_gasket_pentakun}, has \s-s structures. For more details on these structures, see Examples~\ref{eg:sierp_gasket} and \ref{eg:pentakun}.

It is shown in \cite{Kig01} (Proposition~1.3.3) that if $\sss = (K, \alb, \set{F_i}_{i \in \alb})$ is a \s-s structure on $K$, then $\pi$ is unique. Therefore, given a \s-s structure on $K$, we can discuss its surjection $\pi$.

We are interested in the \s-s structures on the limit space of a contracting action.

\begin{cond}
\label{cond:self-similar_eq_rel}
A continuous surjection $\pi: \alb^{-\omega} \to K$ is said to satisfy this condition if $\pi (\ldots x_2 x_1) = \pi (\ldots y_2 y_1)$ implies that $\pi (\ldots x_2 x_1 i) = \pi (\ldots y_2 y_1 i)$ for each $i \in \alb$, and consequently $\pi (\ldots x_2 x_1 w) = \pi (\ldots y_2 y_1 w)$ for each $w \in \alb^*$.
\end{cond}

\begin{prop}
\label{prop:self-similar_eq_rel}
Let $\lims$ be the limit space of a contracting action $(G, \alb)$. There exists a \s-s structure $\sss = (\lims, \alb, \set{F_i}_{i \in \alb})$ on $\lims$, such that the associated continuous surjection $\pi$ is the quotient map $p: \alb^{-\omega} \to \lims$, if and only if the quotient map $p$ satisfies Condition~\ref{cond:self-similar_eq_rel}.
\end{prop}

\begin{proof}
Suppose first that $p$ satisfies Condition~\ref{cond:self-similar_eq_rel}; then $p (\ldots x_2 x_1) = p (\ldots y_2 y_1)$ implies that $p (\ldots x_2 x_1 i) = p (\ldots y_2 y_1 i)$ for each $i \in \alb$. As a consequence, if we define $F_i$ by $F_i \circ p = p \circ i$, we see that $F_i$ is well-defined and continuous. Moreover, $F_i$ is injective for each $i$. Indeed, if $F_i (p (\ldots x_2 x_1) ) = F_i (p (\ldots y_2 y_1) )$, then $p (\ldots x_2 x_1 i) = p (\ldots y_2 y_1 i)$, and so by Proposition~\ref{prop:lim_sp_eq_rel}, $p (\ldots x_2 x_1) = p (\ldots y_2 y_1)$. Therefore, $\sss = (\lims, \alb, \set{F_i}_{i \in \alb})$ is a \s-s structure on $\lims$ with $p$ being the associated continuous surjection.

Conversely, suppose that there exists a \s-s structure $\sss = (\lims, \alb, \set{F_i}_{i \in \alb})$ on $\lims$, such that the associated continuous surjection is $p$. Suppose that $p (\ldots x_2 x_1) = p (\ldots y_2 y_1)$. Since $F_w$ is well-defined for every $w \in \alb^*$, we see that
\[
p (\ldots x_2 x_1 w) = F_w (p (\ldots x_2 x_1) ) = F_w (p (\ldots y_2 y_1) ) = p (\ldots y_2 y_1 w)\text{,}
\]
and so Condition~\ref{cond:self-similar_eq_rel} is satisfied.
\end{proof}

Hereafter, if we say that a limit space $\lims$ has a \s-s structure, we mean that there is a \s-s structure on $\lims$ such that the associated continuous surjection $\pi$ is $p$. In particular, we shall refer to the \s-s structure defined in the proof above as \emph{the} \s-s structure on the limit space $\lims$.

We now wish to investigate which contracting actions have a limit space with a \s-s. We shall see that they are exactly those contracting actions satisfying the following condition.

\begin{cond}
\label{cond:quasi-monocarpic}
A contracting action $(G, \alb)$ is said to satisfy this condition if for every left-infinite path $e = \ldots e_2 e_1$ in the nucleus ending at a non-trivial state and for every $w \in \alb^*$, there exists a left-infinite path $f = \ldots f_2 f_1$ in the nucleus ending at a state $g$, such that the label of the edge $f_n$ is the same as the label of $e_n$, and $g (w) = w$.
\end{cond}

\begin{thm}
\label{thm:quasi-monocarpic}
The limit space $\lims$ of a contracting action $(G, \alb)$ has a \s-s structure if and only if $(G, \alb)$ satisfies Condition~\ref{cond:quasi-monocarpic}.
\end{thm}

\begin{proof}
We are to show that $(G, \alb)$ satisfies Condition~\ref{cond:quasi-monocarpic} if and only if the quotient map $p: \alb^{-\omega} \to \lims$ satisfies Condition~\ref{cond:self-similar_eq_rel}; then we can apply Proposition~\ref{prop:self-similar_eq_rel} to arrive at the desired conclusion.

Suppose first that $(G, \alb)$ satisfies Condition~\ref{cond:quasi-monocarpic}. Let $p (\ldots x_2 x_1) = p (\ldots y_2 y_1)$, then $\ldots x_2 x_1$ is asymptotically equivalent to $\ldots y_2 y_1$. Therefore, there exists a left-infinite path $\ldots e_2 e_1$ within the nucleus passing through the states $\ldots g_2 g_1 g_0$, where the label of the edge $e_n$ is $(x_n, y_n)$.

If $g_0 = 1$, then it is evident that $\ldots x_2 x_1 w$ is asymptotically equivalent to $\ldots y_2 y_1 w$ for all $w \in \alb^*$, and so $p (\ldots x_2 x_1 w) = p (\ldots y_2 y_1 w)$.

If $g_0 \neq 1$, then by Condition~\ref{cond:quasi-monocarpic}, for each $w \in \alb^*$, there exist a state $h \in \nucl$ and a left-infinite path $\ldots f_2 f_1$ within the nucleus ending at $h$, such that the label of $f_n$ is also $(x_n, y_n)$, and that $h (w) = w$. Then $\ldots x_2 x_1 w$ is asymptotically equivalent to $\ldots y_2 y_1 w$, and so $p (\ldots x_2 x_1 w) = p (\ldots y_2 y_1 w)$. Therefore, $p$ satisfies Condition~\ref{cond:self-similar_eq_rel}.

Conversely, suppose that $p$ satisfies Condition~\ref{cond:self-similar_eq_rel}. Let $e = \ldots e_2 e_1$ be a left-infinite path in the nucleus ending at a non-trivial state, where the label of the edge $e_n$ is $(x_n, y_n)$. Then $p (\ldots x_2 x_1) = p (\ldots y_2 y_1)$. By Condition~\ref{cond:self-similar_eq_rel}, $p (\ldots x_2 x_1 w) = p (\ldots y_2 y_1 w)$ for every $w \in \alb^*$. Thus, there exists a left-infinite path $f = \ldots f_2 f_1$ in the nucleus ending at the state $g$, such that the label of the edge $f_n$ is also $(x_n, y_n)$, and $g (w) = g$.
\end{proof}

Henceforth, whenever we discuss a limit space with a \s-s structure, we shall use $\pi$ to denote both the surjection and the quotient map.

For a limit space with a \s-s structure, the notion of a cell and a tile coincides.

\begin{prop}
\label{prop:cell_tile}
Let $\sss = (\lims, \alb, \set{F_i}_{i \in \alb})$ be the \s-s structure on the limit space $\lims$ of a contracting action $(G, \alb)$ satisfying Condition~\ref{cond:quasi-monocarpic}. Then for each $w \in \alb^*$, the $n^{\text{th}}$ level tile $\tile_w$ and the $n^{\text{th}}$ level cell $F_w (\lims)$ are the same set.
\end{prop}

We end this section by proving a strengthened version of Proposition~4.4 of \cite{NT08}.

\begin{prop}
Let $\sss = (\lims, \alb, \set{F_i}_{i \in \alb})$ be the \s-s structure on the limit space $\lims$ of a contracting action $(G, \alb)$ satisfying Condition~\ref{cond:quasi-monocarpic}. Then the restriction of the induced shift map $\shift: \lims \to \lims$ onto the tile $\tile_{wi} = F_{wi} (\lims)$ is equivalent to $F_i^{-1}: \tile_{wi} \to \tile_{w}$ for every $w \in \alb^*$ and $i \in \alb$.

Therefore, the restriction of the induced shift map $\shift$ onto the tile $\tile_{wi}$ is a homeomorphism $\shift: \tile_{wi} \to \tile_{w}$ for every $w \in \alb^*$ and $i \in \alb$. In particular, the tiles are homeomorphic to the limit space $\lims$.
\end{prop}

\begin{proof}
It is evident that the restriction of $\shift: \lims \to \lims$ onto the tile $\tile_{wi}$ is the map $\shift: \tile_{wi} \to \tile_{w}$. Consider now the restriction of $F_i$ onto $\tile_{w}$, which is the map $F_i: \tile_{w} \to \tile_{wi}$. On the set $\alb^{-\omega} w i$ of left-infinite words that end in $w i$, we have that $F_i \circ \shift \circ \pi = F_i \circ \pi \circ \sigma = \pi \circ i \circ \sigma = \pi$, and so we see that $\shift$ and $F_i$ are inverses of each other.

The map $F_i$ is continuous by definition. Since $\shift$ is bijective and continuous, and its inverse is also continuous, we obtain that $\shift: \tile_{wi} \to \tile_{w}$ is a homeomorphism.
\end{proof}

\section{The limit space of a \pcf\ action}
\label{sec:limit_space_pcf_action}

We now turn to a class of self-similar structures that is important in analysis on fractals, namely post-critically finite (or \pcf) structures. The criterion for the limit space of a contracting action to be \fr\ has been shown in \cite{BN03}, and the main result in this section (Theorem~\ref{thm:pcf_pcf}) is to apply this result to the case of limit spaces with a \s-s structure. We first follow \cite{BN03} and adopt the following definition:

\begin{defn}[\cite{BN03} Definition~5.1]
A contracting action is said to be \emph{post-critically finite}, or \pcf\ for short, if there exists only a finite number of left-infinite paths in the Moore diagram of its nucleus which end at a non-trivial state.
\end{defn}

We then follow \cite{Kig93, Kig01} and make the following definitions:

\begin{defn}[\cite{Kig93} Definition~1.5, \cite{Kig01} Definition~1.3.4]
Let $\sss=(K, \alb, \set{F_i}_{i\in \alb})$ be a \s-s structure on $K$. The \emph{critical set} $\crits$ of $\sss$ is defined by $\crits(\sss) = \linebreak \pi^{-1} (\bigcup_{i, j\in \alb, i \neq j} (K_i \cap K_j) )$, and the \emph{post-critical set} $\pcrits$ is defined by $\pcrits(\sss) = \bigcup_{n\geq1} \sigma^n (C)$, where $\sigma$ is the shift operator on $\alb^{-\omega}$.
\end{defn}

\begin{defn}[\cite{Kig93} Definition~1.12, \cite{Kig01} Definition~1.3.13]
A \s-s structure is said to be \emph{post-critically finite}, or \pcf\ for short, if its post critical set $\pcrits$ is finite.
\end{defn}

To prove our main result of this section, we use a lemma from \cite{BN03}. We need the notion of a finitely ramified set to understand the lemma.

\begin{defn}
The limit space of a contracting action is said to be \emph{\fr\ in the group-theoretical sense}, or simply \emph{\fr}, if the intersection of every two distinct tiles of the same level is finite.

A \s-s structure is said to be \emph{\fr\ in the fractal sense}, or simply \emph{\fr}, if the intersection of every two distinct cells of the same level is finite.
\end{defn}

It is a standard result that a \pcf\ \s-s structure is \fr. A \fr\ limit space (in the group-theoretical sense) is what \cite{BN03} calls a \emph{\pcf\ limit space}. As we shall see later, it is true that the \s-s structure of a limit space is \pcf\ if and only if the limit space is \fr; however, to avoid confusion with the notion of a \pcf\ fractal (in the \cite{Kig93, Kig01} sense), we shall not use the terminology introduced by \cite{BN03}.

By Proposition~\ref{prop:cell_tile}, the \s-s structure $\sss=(\lims, \alb, \set{F_i}_{i\in \alb})$ of the limit space of a contracting action is \fr\ (in the fractal sense) if and only if $\lims$ is \fr\ (in the group-theoretical sense).

We now quote the lemma from \cite{BN03}.

\begin{lem}[\cite{BN03} Corollary~4.2]
\label{lem:pcf group_fin_ram_lim_sp}
The limit space $\lims$ of a contracting action $(G, \alb)$ is \fr\ if and only if $(G, \alb)$ is \pcfp
\end{lem}

Our result justifies the use of the terminology ``\pcf''\ in the ``\pcf\ action'' in Lemma~\ref{lem:pcf group_fin_ram_lim_sp}. Before we state our main result, we first prove a useful proposition, the proof of which is used in \cite{BN03} Corollary~4.2.

\begin{prop}
\label{prop:no_inf_iden}
Let $\lims$ be the limit space of a contracting action $(G, \alb)$. Then for every point $a \in \lims$, the set $p^{-1} (a)$ is finite.
\end{prop}

\begin{proof}
We prove that each asymptotic equivalence class of a contracting action has at most $\abs{\nucl}$ elements, and so the quotient map $p$ cannot map infinitely many elements to a point in $K$.

Given a sequence $\ldots x_2 x_1 \in \alb^{-\omega}$, we denote by $E$ the set of all left-infinite paths $\ldots e_2 e_1$ passing through the states $\ldots g_2 g_1 g_0$ within the nucleus, where the label of the edge $e_n$ is $(x_n, y_n)$ for some $y_n \in \alb$. We know that since $g_{n-1}=g_n |_{x_n}$, the state $g_{n-1}$ in the path is uniquely determined by $g_n$ and $x_n$. This shows that given two distinct paths $e =\ldots e_2 e_1$ and $f = \ldots f_2 f_1$ in $E$, if $e_k = f_k$ for some $k$, then $e_n = f_n$ for all $n \leq k$. Consequently, there exists a positive integer $N_{ef}$ such that $e_n \neq f_n$ for all $n \geq N_{ef}$.

Suppose there exist more than $\abs{\nucl}$ distinct left-infinite paths in $E$; then we can choose a set $F$ of $\abs{\nucl} + 1$ distinct paths in $E$. Let $N=\max_{e, f \in F} N_{ef}$, then for every pair of $e, f \in F$, $e_n \neq f_n$ for all $n \geq N$. This is a contradiction since it implies that there are more than $\abs{\nucl}$ states in the nucleus.
\end{proof}

\begin{thm}
\label{thm:pcf_pcf}
The \s-s structure $\sss=(\lims, \alb, \set{F_i}_{i\in \alb})$ on the limit space $\lims$ of a contracting action $(G, \alb)$ is \pcf\ if and only if $(G, \alb)$ is \pcfp
\end{thm}

\begin{proof}
Suppose first that $\sss$ is \pcfp Then $\sss$ is \fr, and so $\lims$ is finitely ramified. Therefore, by Lemma~\ref{lem:pcf group_fin_ram_lim_sp}, $(G, \alb)$ is \pcf

Conversely, suppose that $\sss$ is not \pcfp We can assume that $\sss$ is \fr, for otherwise $\lims$ is not \fr, and we can apply Lemma~\ref{lem:pcf group_fin_ram_lim_sp} to obtain that $(G, \alb)$ is not \pcfp In particular, since the image $\pi (\crits)$ of the critical set is the intersection of cells of the first level, we see that $\pi (\crits)$ is finite.

Moreover, if the critical set $\crits$ is infinite, then there exists at least one point $a \in \pi (\crits)$ such that the set $\pi^{-1} (a)$ is infinite. This is impossible, since by Proposition~\ref{prop:no_inf_iden}, $\lims$ cannot be a limit space. Therefore, $\crits$ is finite.

Now since $\pcrits$ is infinite, there exists at least one element $x \in \crits$ such that the shift map $\sigma$ generates infinitely many distinct elements of $\alb^{-\omega}$ when repeatedly applied to $x$. Since $x \in \crits$, there exists some $y \in \crits$ such that $x_1 \neq y_1$ and $x$ and $y$ are asymptotically equivalent; thus there exists a left-infinite path $\ldots e_2 e_1$ passing through the states $\ldots g_2 g_1 g_0$ within the nucleus, where the label of the edge $e_n$ is $(x_n, y_n)$ and $g_n$ is non-trivial for all $n \geq 1$.

We now show that if $i \neq j$, then the left-infinite paths $\ldots e_{i+1} e_i$ and $\ldots e_{j+1} e_j$, ending at the states $g_{i-1}$ and $g_{j-1}$ respectively, are two distinct paths. Without loss of generality, we can assume $i > j$; if the two paths are identical, then we have that $e_{m+(i-j)} = e_m$ for all $m \geq j$. This would imply that $x_{m+(i-j)}=x_m$ for all $m \geq j$, and consequently $\sigma^{m+(i-j)} (x) = \sigma^m (x)$ for all $m \geq j-1$. But this means that the shift map $\sigma$ only generates at most (including $x$ itself) $(j-1)+(i-j)=i-1$ distinct elements of $\alb^{-\omega}$, contradicting the definition of $x$. Therefore the two paths are distinct.

Also, since for $n \geq 1$ each $g_n$ is a non-trivial state in the nucleus, which is finite because the action is contracting, it follows that there exists an infinite sequence $\set{n_k}$ such that $g_{n_k} = g$ for some non-trivial state $g$ in the nucleus and for all $k$. If we now consider the left-infinite paths $\ldots e_{n_k+2} e_{n_k+1}$, we see that these paths are pairwise distinct and all end at the state $g$, and so there exist infinitely many left-infinite paths in the Moore diagram of the nucleus which end at a non-trivial state.
\end{proof}

\begin{cor}
The \s-s structure $\sss=(\lims, \alb, \set{F_i}_{i\in \alb})$ on the limit space $\lims$ of a contracting action $(G, \alb)$ is \pcf\ if and only if it is \fr. In other words, a \s-s structure that is \fr\ but not \pcf\ cannot be a \s-s structure on the limit space of a contracting action.
\end{cor}

We already know a certain class of fractals that are \fr\ but not \pcf\, that cannot arise as the limit space of any contracting action. In particular, being \fr\ implies that the image $\pi (\crits)$ of the critical set is finite, while not being \pcf\ implies that the post-critical set $\pcrits$ is infinite. If this is the case, Proposition~\ref{prop:no_inf_iden} implies that fractals with either infinite $\crits$ or finite $\pi (\pcrits)$ cannot be the result of a limit space. A simple example of such a fractal is the diamond fractal, which has been discussed in \cite{BCD+08, HK10}.

However, we now have a new class of fractals that are not limit spaces of contracting actions, namely those \s-s sets that are \fr\ and satisfy that
\begin{enumerate}
\item $\crits$ is finite but $\pi (\crits)$ is infinite; and
\item $\pcrits$ is infinite but $\pi (\pcrits)$ is finite.
\end{enumerate}
An example of such fractals is the Kameyama fractal, introduced in \cite{Kam00} in a different setting and discussed in \cite{Hve05}.

Combining these results with those in the last section, we also have the following result.

\begin{cor}
\label{cor:cond_pcf_lim_sp}
The limit space $\lims$ of a contracting action $(G, \alb)$ has a \pcf\ \s-s structure if and only if $(G, \alb)$ satisfies Condition~\ref{cond:quasi-monocarpic} and is \pcfp
\end{cor}

Finally, we look at the related notion of a strictly \pcf\ group, which is defined and discussed in \cite{NT08}. Its definition requires the notion of bounded automata, which is first introduced in \cite{Sid00}. We show that the limit space of a finitely generated strictly \pcf\ group indeed has a \pcf \s-s structure.

\begin{defn}
An automorphism $g \in \Aut \alb^*$ is said to be \emph{bounded} if the Moore diagram of the set $\set{g|_w: w \in \alb^*}$ is finite and its oriented cycles consisting of non-trivial elements are disjoint and not connected by directed paths.
\end{defn}

\begin{defn}[\cite{NT08} Definition~4.2]
A \s-s group $(G, \alb)$ is said to be \emph{strictly \pcf} if and only if it is a subgroup of the group $\C{B} (\alb)$  of bounded automorphisms and every element of the nucleus of $G$ changes at most one letter in every word $w \in \alb^*$.
\end{defn}

The fact that the set $\C{B} (\alb)$ of bounded automorphisms is indeed a group follows from the following theorem in \cite{BN03, Nek05}, which also shows the relationship between $\C{B} (\alb)$ and \pcf\ groups.

\begin{thm}[\cite{BN03} Theorem~5.3, \cite{Nek05} Corollary~3.9.8]
\label{thm:fin_gen_bounded}
The set $\C{B} (\alb)$ of all bounded automorphisms of the tree $\alb^*$ is a group.

A finitely generated \s-s automorphism group $G$ of the tree $\alb^*$ is a \pcf\ group if and only if it is a subgroup of $\C{B} (\alb)$. In particular, every finitely generated \s-s subgroup of $\C{B} (\alb)$ is contracting.
\end{thm}

\begin{cor}
The \s-s structure on the limit space $\lims$ of a finitely generated strictly \pcf\ group $(G, \alb)$ is \pcfp
\end{cor}

\begin{proof}
If $(G, \alb)$ is strictly \pcf, then it is a subgroup of $\C{B} (\alb)$. By Theorem~\ref{thm:fin_gen_bounded}, it is \pcfp
Therefore, by Corollary~\ref{cor:cond_pcf_lim_sp}, we need only to show that if $(G, \alb)$ is strictly \pcf, then it satisfies Condition~\ref{cond:quasi-monocarpic}. This follows from the assumption that every element of the nucleus of $G$ changes at most one letter in every word $w \in \alb^*$. Indeed, for every $g \in \nucl$ and $w \in \alb^*$ such that $g (w) \neq w$, it follows that $g|_w (v) = v$ for all $v \in \alb^*$, or in other words $g|_w = 1$. This implies that if $\ldots e_2 e_1$ is any left-infinite path passing through the states $\ldots g_2 g_1 g_0$ such that the label of some $e_n$ is $(x_n, y_n)$, where $x_n \neq y_n$, then $g_0 = 1$. (In fact, we always have that $n=1$.) Consequently, there are no left-infinite paths with non-trivial labels that do not end at the trivial state, and so Condition~\ref{cond:quasi-monocarpic} is trivially satisfied.
\end{proof}

For an example of a \pcf\ action that satisfies Condition~\ref{cond:quasi-monocarpic} but is not strictly \pcf, see Example~\ref{eg:nonstrict}.

\section{From a \pcf\ \s-s structure}
\label{sec:construction}

In this section, we shall be concerned with the construction of a contracting action whose limit space has a given \pcf\ \s-s structure. More precisely, we shall construct a contracting action whose limit space has a \s-s structure that is isomorphic to a given \pcf\ \s-s structure in the following sense:

\begin{defn}[\cite{Kig01} Definition~1.3.2]
Let $\sss_j = (K_j, \alb_j, \{F_i^{(j)}\}_{i \in \alb_j} )$ be \s-s structures for $j = 1,2$, and let $\pi_j: \alb_j^{-\omega} \to K_j$ be the associated continuous surjections. The \s-s structures $\sss_1$ and $\sss_2$ are said to be \emph{isomorphic} if there exists a bijective map $\rho: \alb_1 \to \alb_2$ such that $\pi_2 \circ \iota_\rho \circ \pi_1^{-1}$ is a well-defined homeomorphism between $K_1$ and $K_2$, where $\iota_\rho: \alb_1^{-\omega} \to \alb_2^{-\omega}$ is the natural bijective map defined by $\iota(\ldots x_2 x_1) = \ldots \rho(x_2) \rho (x_1)$.
\end{defn}

Notice that given a set with a \s-s structure in general, it is possible that it is not the limit space of any contracting action. In particular, Proposition~\ref{prop:lim_sp_eq_rel} implies that for any construction to have a hope of success, the surjection $\pi$ of the \s-s structure must be such that if $\pi (\ldots x_2 x_1) = \pi (\ldots y_2 y_1)$, then $\pi (\ldots x_{n+1} x_n) = \pi (\ldots y_{n+1} y_n)$. An equivalent requirement is that the induced shift map $\shift: \lims \to \lims$ defined by $\shift \circ \pi = \pi \circ \sigma$ is well-defined. For example, the usual \s-s structure on the Sierpi\'nski gasket does not satisfy this requirement, although another \s-s structure on it does; see Example~\ref{eg:sierp_gasket} for a more detailed description.

At the same time, a property of \s-s structures will also be useful. A close inspection of the proof of Proposition~\ref{prop:self-similar_eq_rel} reveals that for any \s-s structure, the surjection $\pi$ must satisfy Condition~\ref{cond:self-similar_eq_rel}. This leads us to the following lemma.

\begin{lem}\label{lem:shift_conditions}
Let $\sss = (K, \alb, \{F_i\}_{i \in \alb} )$ be a \s-s structure on the limit space of a contracting action $(G, \alb)$. Then the associated surjection $\pi: \alb^{-\omega} \to K$ satisfies the following conditions:
\begin{enumerate}
\item $\pi (\ldots x_2 x_1) = \pi (\ldots y_2 y_1)$ implies that $\pi (\ldots x_{n+1} x_n) = \pi (\ldots y_{n+1} y_n)$, for all $n \in \B{Z}^+$; and
\item $\pi (\ldots x_2 x_1) = \pi (\ldots y_2 y_1)$ implies that $\pi (\ldots x_2 x_1 w) = \pi (\ldots y_2 y_1 w)$, for all $w \in \alb^*$.
\end{enumerate}
\end{lem}

In this section, we restrict ourselves to considering only \s-s structures that are isomorphic to \s-s structures on limit spaces. In particular, for the rest of this section, we shall take for granted the existence of the shift map $\shift$, and that the associated surjection $\pi$ satisfies the conditions in Lemma~\ref{lem:shift_conditions}.

Given a \pcf\ \s-s structure $\sss = (K, \alb, \set{F_i}_{i \in \alb} )$, the surjection $\pi$ defines an equivalence relation on $\alb^{-\omega}$. We now describe a scheme to systematically write down the equivalence classes induced by $\pi$.

Suppose $\pi (\ldots x_2 x_1) = \pi (\ldots y_2 y_1)$, and $x_k = y_k$ for all $k < N$ and $x_N \neq y_N$. By Condition~(1) above, we have that $\pi (\ldots x_{N+1} x_N) = \pi (\ldots y_{N+1} y_N)$. Then by Condition~(2), $\pi (\ldots x_{N+1} x_N w) = \pi (\ldots y_{N+1} y_N w)$ for all $w \in \alb^*$, which accounts for the fact that $\pi (\ldots x_2 x_1) = \pi (\ldots y_2 y_1)$. The equivalence relation can be completely characterized by equations of the form
\[
\pi (\ldots x_2 x_1 w) = \pi (\ldots y_2 y_1 w)\text{,}
\]
where $x_1 \neq y_1$. Moreover, the equation above implies that $\ldots x_2 x_1, \ldots y_2 y_1 \in \crits$. Since $\sss$ is \pcf, $\crits$ must be finite. Consequently, the equivalence relation can be characterized by finitely many such equations.

The fact that $\sss$ is \pcf\ also implies that elements in $\crits$ have a recurring tail. If we denote $\recur{z} = \ldots z z z$ where $z = z_k \ldots z_2 z_1 \in \alb^*$, then the equivalence relation can be characterized by finitely many equations of the form
\[
\pi (\recur{z} x_n \ldots x_2 x_1 w) = \pi (\recur{z}^\prime y_m \ldots y_2 y_1 w) \text{.}
\]
We shall now show that $m = n$ in the equation above. Otherwise, without loss of generality, we can let $m > n$. Then by Condition~(1), $\pi (\recur{z}) = \pi (\recur{z}^\prime y_m \ldots y_{n+1})$. Then
\[
\pi (\recur{z}^\prime y_m \ldots y_{n+1}) = \pi (\recur{z}) = \pi (\recur{z} z) = \pi (\recur{z}^\prime y_m \ldots y_{n+1} z)\text{,}
\]
by Condition (2); likewise,
\[
\pi (\recur{z}^\prime y_m \ldots y_{n+1}) = \pi (\recur{z}^\prime y_m \ldots y_{n+1} z) = \pi (\recur{z}^\prime y_m \ldots y_{n+1} z z) = \ldots \text{,}
\]
which would imply that $\crits$ is infinite, again a contradiction to the fact that $\sss$ is \pcfp A similar argument shows that we must have $z = z^\prime$. Therefore, the equivalence relation can be characterized by finitely many equations of the form
\[
\pi (\recur{z} x_n \ldots x_2 x_1 w) = \pi (\recur{z} y_n \ldots y_2 y_1 w)\text{.}
\]
Another similar argument shows that $z_k \neq x_n$ and $z_k \neq y_n$. We shall assume that $z$ is the shortest recurring word, so that it is impossible to write $z = v^n$ for any $n > 1$ and $v \in \alb^*$.

Next, we show that if $\pi (\recur{z} x_n \ldots x_2 x_1 w) = \pi (\recur{z} y_n \ldots y_2 y_1 w)$, then $\pi (\recur{z} \xi_n \ldots \xi_2 \xi_1 w) = \pi (\recur{z} x_n \ldots x_2 x_1 w)$ whenever $\xi_j \in \set{x_j, y_j}$ for all $j$. We shall show this by induction on $n$. The base case $n=1$ is trivial. Suppose this is true for $n = m$, then by Condition~(2),
\begin{align*}
\pi (\recur{z} \xi_{m+1} \ldots \xi_2 x_1 w) &= \pi (\recur{z} x_{m+1} \ldots x_2 x_1 w) \\ &= \pi (\recur{z} y_{m+1} \ldots y_2 y_1 w)\\ & = \pi (\recur{z} \xi_{m+1} \ldots \xi_2 y_1 w)\text{,}
\end{align*}
and therefore it is true for $n = m+1$. As a consequence, we have that $\pi (\recur{z} \xi_n \ldots \xi_2 \xi_1 w) = \pi (\recur{z} \zeta_n \ldots \zeta_2 \zeta_1 w)$ whenever $\xi_j, \zeta_j \in \set{x_j, y_j}$ for all $j$. At this point, having classified the equivalence classes induced by the associated surjection $\pi$ of a \pcf\ \s-s structure, we can finally write the equivalence classes in the form
\[
\set{\recur{z} \zeta_n \ldots \zeta_2 \zeta_1 w \: | \: z, w \in \alb^*, \:\zeta_j \in S_j}
\]
for fixed $w, z \in \alb^*$, and some collection of sets $S_j \subset \alb$. We introduce the shorthand
\[
\recur{z} S_n \ldots S_2 S_1 w \text{,}
\]
to represent the equivalence class above.

Up to left shifts, all the equivalence classes are determined by those of the form $\recur{z} S_n \ldots S_1$, where $S_1$ contains more than one element. Notice that if $\alpha$ is in the image $\pi(\crits)$ of the critical set, then $\pi^{-1}(\alpha) = \recur{z}S_n \ldots S_1$, for some $z\in X^*$ and $S_1$ with more than one element, since $\alpha$ is in the union of the intersections of the first level cells of the space $K$; conversely, the image of every equivalence class of the form $\recur{z} S_n \ldots S_1$ under $\pi$ is a single point in $\pi (\crits)$. Therefore, the equivalence classes can be labeled by the finitely many elements of $\pi (\crits)$.

Every equivalence class of the form $\recur{z} S_n \ldots S_2 S_1$ has to satisfy three properties. First, by what we discussed above, we see that if $z = z_k \ldots z_2 z_1$, then $z_k \notin S_n$. Second, if $\recur{z} S_n \ldots S_2 S_1$ is in the list, then by Condition~(1), we must have that $\recur{z} S_n \ldots S_{m+1} S_m$ is also in the list for all $m \leq n$. The third property is the proposition below:

\begin{prop}
Let $\alpha$ and $\beta$ be distinct elements of $\pi(\crits)$ of a \pcf \ \s-s structure, such that $\pi^{-1} (\alpha)$ and $\pi^{-1} (\beta)$ have the same recurring tail $\recur{z}$. Let $\pi^{-1}(\alpha) = \recur{z}S_n \ldots S_2S_1$ and $\pi^{-1}(\beta) = \recur{z}T_m\ldots T_2 T_1$. If $S_{n-k} = T_{m-k}$ for all $0 \leq  k < N$, then either $S_{n-N} = T_{m-N}$ or $S_{n-N} \cap T_{m-N} = \varnothing$.
\end{prop}

\begin{proof}
Suppose $S_{n-N} \cap T_{m-N} \neq \varnothing$, and let $x \in S_{n-N} \cap T_{m-N}$. Then for all $s_{n-N} \in S_{n-N}$ and $t_{m-N} \in T_{m-N}$, we have 
\begin{align*}
\pi (\recur{z} \xi_n \ldots \xi_{n-N+1} s_{n-N} \xi_{n-N-1} \ldots \xi_2 \xi_1) & = \pi (\recur{z} \xi_n \ldots \xi_{n-N+1} x \xi_{n-N-1} \ldots \xi_2 \xi_1)\\
& = \pi (\recur{z} \xi_n \ldots \xi_{n-N+1} t_{m-N} \xi_{n-N-1} \ldots \xi_2 \xi_1) \text{,}
\end{align*}
whenever $\xi_j \in S_j$. Therefore, we see that $T_{m-N} \subset S_{n-M}$. Similarly, we obtain that $S_{n-N} \subset T_{m-N}$, and so $S_{n-N} = T_{m-N}$.
\end{proof}

We can now construct the desired contracting group. For each $\alpha \in \pi (\crits)$, we can write $\pi^{-1} (\alpha) = \recur{z} S_n \ldots S_2 S_1 $, where $z=z_k \ldots z_2 z_1$. If $S_j = \{ s_1^{(j)}, s_2^{(j)}, \ldots, s_m^{(j)} \}$ with $s_i^{(j)} < s_{i+1}^{(j)}$, we define $\sigma_{S_j}$ to be the permutation $(s_1^{(j)} s_2^{(j)} \ldots s_m^{(j)})$. We define $n+k-1$ group elements as follows:

We define $g_{\alpha,1}$ by the wreath recursion $\sigma_{S_1}(1,\ldots,1)$. For $2 \leq j \leq n-1$, we define $g_{\alpha, j}$ to be the element whose action on $x \in \alb$ is given by
\[
g_{\alpha, j} \cdot x =
\begin{cases}
\sigma_{S_j}(x) \cdot g_{\alpha, j-1} & \quad \text{if } x \in S_j \text{,}\\
x \cdot 1 & \quad \text{if } x \notin S_j \text{,}
\end{cases}
\]
so that we have the wreath recursion $g_{\alpha,j} = \sigma_{S_j} (1, \ldots, g_{\alpha, j-1}, \ldots, g_{\alpha, j-1}, \ldots, 1)$.

For $j = n$, we define $g_{\alpha, n}$ by
\[
g_{\alpha, n} \cdot x =
\begin{cases}
\sigma_{S_n}(x)\cdot g_{\alpha, n-1} & \quad \text{if } x \in S_n \text{,} \\
x \cdot g_{\alpha, n+k-1} & \quad \text{if } x = z_k \text{,} \\
x \cdot 1 & \quad \text{otherwise,}
\end{cases}
\]
so that we have the wreath recursion $g_{\alpha, n} = \sigma_{S_n} (1, \ldots, g_{\alpha, n-1}, \ldots, g_{\alpha, n-1}, \ldots, g_{\alpha, n+k-1}, \ldots, 1)$.

Finally, for $n+1 \leq j \leq n+k-1$, we define $g_{\alpha, j}$ by
\[
g_{\alpha, j} \cdot x =
\begin{cases}
x \cdot g_{\alpha, j-1} & \quad \text{if } x = z_{j-n} \text{,}\\
x \cdot 1 & \quad \text{if } x \neq z_{j-n} \text{,}
\end{cases}
\]
so that we have the wreath recursion $g_{\alpha, j} = (1, \ldots, g_{\alpha, j-1}, \ldots, 1)$.

We call $g_{\alpha, 1}, \ldots, g_{\alpha, n-1}$ \emph{Type~I generators} and $g_{\alpha, n} \ldots, g_{\alpha, n+k-1}$ \emph{Type~II generators}.

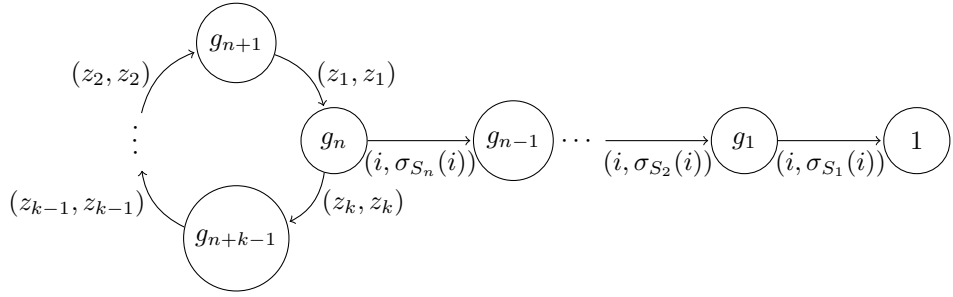
\begin{figure}
\begin{tikzpicture}[->, shorten >=1pt, node distance=1.4cm]
\node[state] at (0:1.3)      (gn)     {$g_n$};
\node[state] at (90:1.3)     (gn+1)   {$g_{n+1}$};
\node[state] at (270:1.3)    (gn+k-1) {$g_{n+k-1}$};
\node[state] [right=of gn]   (gn-1)   {$g_{n-1}$};
\node        at (4.55,0)      (a)      {$\cdots$};
\node[state] [right=of a]    (g1)     {$g_1$};
\node[state] [right=of g1]   (1)      {$1$};
\node        at (173:1.35)   (b)      {$\cdot$};
\node        at (180:1.35)   (c)      {$\cdot$};
\node        at (187:1.35)   (d)      {$\cdot$};
\path
(g1)     edge                node [below]           {$(i,\sigma_{S_1}(i))$}     (1)
(a)      edge                node [below]           {$(i,\sigma_{S_2}(i))$}     (g1)
(gn)     edge                node [below]           {$(i,\sigma_{S_n}(i))$}     (gn-1)
         edge [bend left=28] node [right]           {$(z_k,z_k)$}               (gn+k-1)
(gn+1)   edge [bend left=28] node [right]           {$(z_1,z_1)$}               (gn)          
(b)      edge [bend left=28] node [left]            {$(z_2,z_2)$}               (gn+1)
(gn+k-1) edge [bend left=28] node [left]            {$(z_{k-1},z_{k-1})$}       (d);
\end{tikzpicture}
\caption{The generators corresponding to $\alpha = \pi (\recur{z} S_n \ldots S_2 S_1)$}
\label{fig:construction}
\end{figure}

The Moore diagram of these generators, which correspond to a single $\alpha \in \pi (\crits)$, is shown in Figure~\ref{fig:construction}, in which the subscript $\alpha$ has been suppressed. The label $(i, \sigma_{S_j} (i))$ for the arrows out of $g_n$ and the Type~II generators applies to all $i \in S_j$ and only those letters. For each generator, we have suppressed the arrows into the identity element, whose labels are $(i, i)$ for each $i$ that has not been shown in the diagram.

The desired group $G_\sss$ is the group generated by all the elements defined above for all $\alpha \in \pi (\crits)$.

\begin{prop}
\label{prop:construction_pcf}
Let $\sss = (K, \alb, \set{F_i}_{i \in \alb})$ be a \pcf\ \s-s structure. The group $G_\sss$ that is constructed by the method above is a subgroup of $\C{B} (X)$. In particular, $G_\sss$ is contracting and \pcfp
\end{prop}

\begin{proof}
It can easily be seen from the Moore diagram above that all the generators of $G_\sss$ are bounded automorphisms. Since $G_\sss$ is finitely generated, it follows from Theorem~\ref{thm:fin_gen_bounded} that $G_\sss$ is \pcfp 
\end{proof}

For concrete examples illustrating our construction above, see Section~\ref{sec:examples}, and in particular Examples~\ref{eg:alt_grigorchuk}, \ref{eg:pentakun} and \ref{eg:hexakun_snowflake}.

We claim that the \s-s structure $\sss^\prime = (\C{J}_{G_\sss}, \alb, \set{F_i^\prime}_{i \in \alb})$ on the limit space of $G_\sss$ is isomorphic to $\sss$ (Theorem~\ref{thm:construction_isom}). We begin the proof with a lemma, where we show that the associated continuous surjection of a \s-s structure is in fact a quotient map.

\begin{lem}
\label{lem:pi_quot_map}
Let $\sss = (K, \alb, \set{F_i}_{i \in \alb})$ be a \s-s structure. The associated continuous surjection $\pi: \alb^{-\omega} \to K$ of $\sss$ is a quotient map.
\end{lem}

\begin{proof}
We shall show that $\pi$ is a closed map. Since $\alb^{-\omega}$ is compact, every closed subset $C$ is compact. Then $\pi (C)$ is a compact subset of $K$, which is metrizable and thus Hausdorff, implying that $\pi (C)$ is closed.
\end{proof}

\begin{lem}
\label{lem:repres_by_gen}
Let $G_\sss$ be the contracting group constructed by the method above. If $h \in G_\sss$ can be written as a product $g_m \ldots g_2 g_1$ of generators of minimal length, and if $h|_v = h$ for some non-empty word $v = v_n \ldots v_2 v_1 \in \alb^*$, then $g_j (v) = v$ and $g_j |_v = g_j$ for all $j$. In particular, $g_j$ is of Type~II for all $j$, and $h(v) = v$.
\end{lem}

\begin{proof}
We denote by $d(g)$ the minimal length of generators needed to represent an element $g \in G_\sss$. In particular, $d(h) = m$. Notice that, since $h|_x = g_m|_{g_{m-1}\cdots g_1(x)} \cdots g_2|_{g_1(x)}g_1|_x$, and $g_j|_y$ is $1$ or another generating element for all $y \in X$, that $d(h|_x) \leq d(h)$. We see $d(h) = d(h|_v) \leq d(h|_{v_n\ldots v_{n-i}}) \leq d(h)$ for all $0 < i \leq n$, so $d(h|_{v_n\cdots v_{n-i}}) = d(h)$ for all such $i$.

We show that $g_j$ is of Type~II for all $j$. If $g_j$ is of Type~I for some $j \leq m$, then there exists $N$ such that $g_j |_w = 1$ whenever $\abs{w} = N$. Since if $g$ is a generator, then $g|_x$ is either a generator or the identity for all $x \in \alb$, it follows that $d(h|_w) < m$ whenever $\abs{w} = N$. This is a contradiction since $h|_v = h$ implies that for each $k$, there exists some $w \in \alb^k$ such that $d(h|_w) = m$.

By the same argument, each representation of
\[
h|_{v_n} = (g_m \ldots g_1) |_{v_n} = g_m |_{g_{m-1} \ldots g_1 (v_n)} \ldots g_1 |_{v_n}
\]
of length $m$ consists only of Type~II generators, which implies that $g_1 |_{v_n}$ is also of Type~II. This is only possible if $g_1 (v_n) = v_n$, since for any generator $g$ and any $x \in \alb$, $g(x) \neq x$ implies that $g|_x$ is either Type I or the identity. Similarly, $g_j (v_n) = v_n$ and $g_j |_{v_n}$ is of Type~II for all $j$, and
\[
h|_{v_n} = g_m |_{g_{m-1} \ldots g_1 (v_n)} \ldots g_1 |_{v_n} = g_m |_{v_n} \ldots g_1 |_{v_n} \text{.}
\]
Inductively, $g_j (v) = v$ and $g_j |_v$ is of Type~II for each $j$. Moreover, $g_j |_v$ and $g_j$ are generators corresponding to the same $\alpha \in \pi (\crits)$. Since there are only finitely many Type~II generators corresponding to $\alpha$, it follows that $g_j |_{v^l} = g_j$ for some minimal $l$.

Suppose $l > 1$. If $g_j$ corresponds to $\alpha \in \pi (\crits)$, then $\pi^{-1} (\alpha) = \recur{z} S_n \ldots S_2 S_1$, where $z = v^l$. This is a contradiction since we assume $z$ to be the shortest recurring word.
\end{proof}

\begin{lem}
\label{lem:p_pi_same_eq}
Let $\sss = (K, \alb, \set{F_i}_{i \in \alb})$ be a \pcf\ \s-s structure, with $\pi: \alb^{-\omega} \to K$ as the associated surjection. Let $G_\sss$ be the contracting group constructed by the method above, and let $\C{J}_{G_\sss}$ be its limit space. Then the quotient map $p: \alb^{-\omega} \to \C{J}_{G_\sss}$ is such that
\[
p(\ldots x_2 x_1) = p(\ldots y_2 y_1)
\]
if and only if
\[
\pi(\ldots x_2 x_1) = \pi(\ldots y_2 y_1) \text{.}
\]
\end{lem}

\begin{proof}
Suppose first that $\pi(\ldots x_2 x_1) = \pi(\ldots y_2 y_1)$, where $\ldots x_2 x_1 \neq \ldots y_2 y_1$. Without loss of generality, we can instead write $\pi (\recur{z} x_n \ldots x_2 x_1 w) = \pi (\recur{z} y_n \ldots y_2 y_1 w)$, with $z = z_k \ldots z_2 z_1$ and $x_1 \neq y_1$. All but finitely many truncations of the left-infinite word $\recur{z} x_n \ldots x_2 x_1 w$ are of the form $z_l \ldots z_2 z_1 z^N x_n \ldots x_2 x_1 w$, where $1 \leq l \leq k$ and $N \geq 0$. We claim that for each word in this form, there exists an element $g$ such that $g (z_l \ldots z_2 z_1 z^N x_n \ldots x_2 x_1 w) = z_l \ldots z_2 z_1 z^N y_n \ldots y_2 y_1 w$. We shall proceed by induction on $n$.


Consider the case when $n = 1$. There exists $\alpha \in \pi(\crits)$ with $\pi^{-1}(\alpha) = \recur{z} S_1$, such that $x_1, y_1 \in S_1$. By construction, the action of $g_{\alpha, l+1}$ on $z_l \ldots z_2 z_1 z^N x_1 w$ is
\begin{align*}
g_{\alpha, l+1} (z_l \ldots z_2 z_1 z^N x_1 w) & = z_l g_{\alpha, l} (z_{l-1} \ldots z_2 z_1 z^N x_1 w)\\
& = z_l z_{l-1} g_{\alpha, l-1} (z_{l-2} \ldots z_2 z_1 z^N x_1 w)\\
& \hspace{2mm} \vdots\\
& = z_l \ldots z_2 z_1 z^N g_{\alpha, 1} (x_1 w)\\
& = z_l \ldots z_2 z_1 z^N \sigma_{S_1} (x_1) w \text{.}
\end{align*}
There exists $a \in \B{N}$ such that $\sigma_{S_1}^a (x_1) = y_1$. Then
\[
g_{\alpha, l+1}^a (z_l \ldots z_2 z_1 z^N x_1 w) =z_l \ldots z_2 z_1 z^N \sigma_{S_1}^a (x_1) w = z_l \ldots z_2 z_1 z^N y_1 w \text{.}
\]
Suppose now that the statement holds for $n = m$, and consider the case when $n = m + 1$. There exists $\beta \in \pi(\crits)$ such that $\pi^{-1}(\beta) = \recur{z} S_{m+1} \ldots S_2 S_1$, such that $x_k, y_k \in S_k$ for each $k$. By construction, the action of $g_{\beta, l+m+1}$ on $z_l \ldots z_2 z_1 z^N x_{m+1} \ldots x_2 x_1 w$ is
\begin{align*}
g_{\beta, l+m+1} (z_l \ldots z_2 z_1 z^N x_{m+1} \ldots x_2 x_1 w) & = z_l g_{\beta, l+m} (z_{l-1} \ldots z_2 z_1 z^N x_{m+1} \ldots x_2 x_1 w)\\
& \hspace{2mm} \vdots\\
& = z_l \ldots z_2 z_1 z^N g_{\beta, m+1} (x_{m+1} \ldots x_2 x_1 w)\\
& = z_l \ldots z_2 z_1 z^N \sigma_{S_{m+1}} (x_{m+1}) \ldots \sigma_{S_1} (x_1) w
\end{align*}
There exists $b \in \B{N}$ such that $\sigma_{S_1}^b (x_1) = y_1$. Then
\[
g_{\beta, l+m+1}^b (z_l \ldots z_1 z^N x_{m+1} \ldots x_1 w) = z_l \ldots z_1 z^N \sigma_{m+1}^b (x_{m+1}) \ldots \sigma_2^b (x_2) y_1 w \text{.}
\]
By the induction hypothesis, there exists an element $h$ that maps this to
\[
z_l \ldots z_1 z^N y_{m+1} \ldots y_2 y_1 w \text{;}
\]
then $hg_{\beta, l+m+1}^b$ is the desired element. 

Notice that in both cases above, the element fulfilling our claim is independent of $N$; therefore, the set of elements fulfilling our claim for each truncation of the left-infinite word $\recur{z} x_n \ldots x_2 x_1 w$ is finite, and thus $p(\recur{z} x_n \ldots x_1 w) = p(\recur{z} y_n \ldots y_1 w)$.

Conversely, suppose $p (\ldots x_2 x_1) = p (\ldots y_2 y_1)$, where $\ldots x_2 x_1 \neq \ldots y_2 y_1$. Then there exists a left-infinite path $\ldots e_2 e_1$ in $\nucl$ passing through the states $\ldots h^{(2)} h^{(1)} h^{(0)}$, such that the label of $e_k$ is $(x_k, y_k)$. We can represent each $h^{(k)}$ as a product $g_{m_k}^{(k)} \ldots g_1^{(k)}$ of generators with minimal length $m_k$. If $g$ is a generator, then $g|_x$ is either another generator or the identity for all $x \in \alb$; therefore, $m_k \leq m_{k+1}$ for all $k$. At the same time, since the nucleus $\nucl$ is finite, $m_k$ is uniformly bounded from above. Therefore, there exists some minimal $M$ such that $m_k = m$ for some constant $m$ for all $k \geq M$.

By Proposition~\ref{prop:construction_pcf}, $G_\sss$ is \pcfp Therefore, the path $\ldots e_2 e_1$ must have a recurring tail; that is, there exists some minimal $N \geq M$ and some minimal $r \geq 1$ such that $e_{k+r} = e_{k}$ for all $k > N$. Then $h^{(k+r)} = h^{(k)}$ for all $k \geq N$. If we write $z = z_r \ldots z_1 = x_{N+r} \ldots x_{N+1}$, then we have
\[
p \left(\recur{z} x_N \ldots x_2 x_1\right) = p \left(\recur{h^{(N)} (z)} y_N \ldots y_2 y_1\right) \text{.}
\]
Now we also have that $h^{(N)} |_z = h^{(N)}$. By Lemma~\ref{lem:repres_by_gen}, it follows that $g_j^{(N)} (z) = z$ and $g_j^{(N)} |_z = g_j^{(N)}$ for all $j \leq m$, and $h^{(N)}(z) = z$.

Then our original equation becomes
\[
p (\recur{z} x_N \ldots x_2 x_1) = p (\recur{z} y_N \ldots y_2 y_1) \text{,}
\]
where $y_n \ldots y_2 y_1 = g_m^{(N)} \ldots g_1^{(N)} (x_n \ldots x_2 x_1)$.

Consider the generator $g_1^{(N)}$. Since $g_1^{(N)} (z) = z$ and $g_1^{(N)} |_z = g_1^{(N)}$, there exists a left-infinite path $\ldots f_2 f_1$ ending at $g_1^{(N)}$, where the label of $f_{rt+s}$ is $(z_s, z_s)$, for $s < r$. Therefore,
\[
p (\recur{z} x_N \ldots x_2 x_1) = p (\recur{z} g_1^{(N)} (x_N \ldots x_2 x_1) \text{.}
\]
Similarly, we obtain,
\[
p (\recur{z} x_N \ldots x_2 x_1) = p (\recur{z} g_1^{(N)} (x_N \ldots x_2 x_1) = \ldots = p (\recur{z} g_m^{(N)} \ldots g_1^{(N)} (x_N \ldots x_2 x_1)) \text{,}
\]
from which we deduce that $g_j^{(N)} \ldots g_1^{(N)} (x_N) \neq z_r$ for all $j$.

We now show that the above equation will continue to hold if we replace $p$ by $\pi$. For each $j$, we can write $g_j^{(N)} = g_{\alpha_j, l_j}$ for some $\alpha_j \in\pi(\crits)$ with $\pi^{-1}(\alpha_j) = \recur{z^{(j)}} S_{j, n_j} \ldots S_{j, 2} S_{j, 1} $, where $l_j \geq n_j$, since $g_j^{(N)}$ is of Type~II by Lemma~\ref{lem:repres_by_gen}. If $l_j \neq n_j$, then $g_j^{(N)} |_x = 1$ for all $x \neq z_r$. Therefore,
\[
g_j^{(N)} (g_{j-1}^{(N)} \ldots g_1^{(N)} (x_N)) = g_{j-1}^{(N)} \ldots g_1^{(N)} (x_N)
\]
and
\[
g_j^{(N)} |_{(g_{j-1}^{(N)} \ldots g_1^{(N)} (x_N)} = 1 \text{.}
\]
Then
\[
g_j^{(N)} (g_{j-1}^{(N)} \ldots g_1^{(N)} (x_N \ldots x_1)) = g_{j-1}^{(N)} \ldots g_1^{(N)} (x_N \ldots x_1) \text{.}
\]
It is then trivially true that
\[
\pi (\recur{z} g_j^{(N)} (g_{j-1}^{(N)} \ldots g_1^{(N)} (x_N \ldots x_1))) = \pi (\recur{z} g_{j-1}^{(N)} \ldots g_1^{(N)} (x_N \ldots x_1)) \text{.}
\]

Suppose now that $l_j = n_j$; then by construction, $z^{(j)}$ is the unique word of minimal length such that $g_j^{(N)} |_{z^{(j)}} = g_j^{(N)}$. Therefore, we see that $z^{(j)} = z$. If we let $c$ be the smallest number such that $g_j^{(N)} |_{g_{j-1}^{(N)} \ldots g_1^{(N)} (x_N \ldots x_c)} \neq 1$, then by construction,
\begin{gather*}
(g_{j-1}^{(N)} \ldots g_1^{(N)}) |_{x_N \ldots x_{k+1}} (x_k) \in S_{j, n_j-N+k} \text{, and}\\
g_j^{(N)} |_{g_{j-1}^{(N)} \ldots g_1^{(N)} (x_N \ldots x_{k+1})} ((g_{j-1}^{(N)} \ldots g_1^{(N)}) |_{x_N \ldots x_{k+1}} (x_k)) \in S_{j, n_j-N+k} \text{,}
\end{gather*}
whenever $c \leq k \leq N$. But $\recur{z} S_{j, n_j} \ldots S_{j, 2} S_{j, 1} w $ is an equivalence class induced by $\pi$, and this implies that $\recur{z} S_{j, n_j} \ldots S_{j, n_j-N+c} w$ is an equivalence class, and so again we have
\[
\pi (\recur{z} g_j^{(N)} (g_{j-1}^{(N)} \ldots g_1^{(N)} (x_N \ldots x_1))) = \pi (\recur{z} g_{j-1}^{(N)} \ldots g_1^{(N)} (x_N \ldots x_1) \text{.}
\]

Therefore, inductively, we have shown that
\[
\pi (\recur{z} x_N \ldots x_2 x_1) = \pi (\recur{z} y_N \ldots y_2 y_1)) \text{,}
\]
which is what we wanted to prove.
\end{proof}

\begin{thm}
\label{thm:construction_isom}
Let $\sss = (K, \alb, \set{F_i}_{i \in \alb})$ be a \pcf\ \s-s structure, and let $G_\sss$ be the contracting group constructed by the method above. Then the limit space $\C{J}_{G_\sss}$ of $G_\sss$ has a \s-s structure $\sss^\prime = (\C{J}_{G_\sss}, \alb, \set{F_i^\prime}_{i \in \alb})$. Moreover, $\sss^\prime$ is isomorphic to $\sss$.
\end{thm}

\begin{proof}
We know that $\pi$ satisfies Condition~\ref{cond:self-similar_eq_rel}. By Lemma~\ref{lem:p_pi_same_eq}, we see that $p$ also satisfies Condition~\ref{cond:self-similar_eq_rel}. Therefore, by Proposition~\ref{prop:self-similar_eq_rel}, $\C{J}_{G_\sss}$ has a \s-s strucutre $\sss^\prime$.

To show that $\sss^\prime$ is isomorphic to $\sss$, consider the map $p \circ \pi^{-1}$. Given an open set $U$ in $\C{J}_{G_\sss}$, $p^{-1} (U)$ is open because $p$ is continuous. Since by Lemma~\ref{lem:pi_quot_map} $\pi$ is a quotient map, and by Lemma~\ref{lem:p_pi_same_eq} $p^{-1} (U)$ is saturated, it follows that $\pi \circ p^{-1} (U)$ is open. Therefore, $p \circ \pi^{-1}$ is continuous. Since $p$ is by definition a quotient map, we obtain that the inverse is also continuous. Thus, $p \circ \pi^{-1}$ is a well-defined homeomorphism between $\C{J}_{G_\sss}$ and $K$.
\end{proof}

In the construction above, we have only assumed that the induced shift map $\shift$ with regard to the \s-s structure is well-defined, which is a necessary condition for any construction to succeed. By Theorem~\ref{thm:construction_isom}, it is also the sufficient condition; therefore, we have the following theorem.

\begin{thm}
\label{thm:construct_iff_s_exists}
Let $\sss = (K, \alb, \set{F_i}_{i \in \alb})$ be a \pcf\ \s-s structure, with an associated continuous surjection $\pi: \alb^{-\omega} \to K$. Then there exists a contracting action $(G, \alb)$ such that its limit space $\lims$ has a \s-s structure $\sss^\prime = (\lims, \alb, \set{F_i^\prime}_{i \in \alb})$ that is isomorphic to $\sss$, if and only if the induced shift map $\shift: K \to K$ defined by $\shift \circ \pi = \pi \circ \sigma$ is well-defined.
\end{thm}

\section{Examples}
\label{sec:examples}

\begin{eg}[The binary adding machine]
\label{eg:adding_machine}
One of the most basic examples of a \s-s action is the binary adding machine. It is the group $G$ generated by the element $a = (01) (1, a)$, acting on the binary tree (i.e.\ $\alb = \set{0,1}$). The action of $a$ can be thought of as adding 1 to the last digit of the (left-handed) binary representation of a real number.

The nucleus of $G$ is $\set{1,a,a^{-1}}$, which is depicted in Figure~\ref{fig:adding_machine}. The asymptotic equivalence is given by $\recur{0} 1 w \sim \recur{1} 0 w$ for all $w \in \alb^*$ and $\recur{0} \sim \recur{1}$. Therefore, the limit space $\lims$ is homeomorphic to the circle $\B{R} / \B{Z}$, where each point on the circle corresponds to its (left-infinite) binary expansion.

The action $(G, \alb)$ does not satisfy Condition~\ref{cond:quasi-monocarpic}, and therefore its limit space does not have the naturally induced \s-s structure. To see this, notice that $\recur{0} \sim \recur{1}$ but $\recur{0} \not\sim \recur{1} 0$, and so $F: v \mapsto v 0$ is not a well defined map on $\alb^{-\omega}$, which shows the non-existence of a \s-s structure on $\B{R} / \B{Z}$. Therefore, although $G$ is clearly \pcf, its limit space does not have a \pcf\ \s-s structure.
\end{eg}

\begin{eg}[Sierpi\'{n}ski gasket]
\label{eg:sierp_gasket}
The Sierpi\'nski gasket is typically defined, e.g.\ in \cite{Kig01}, as the unique non-empty compact space $K \subset \B{C}$ that is invariant under the injections $f_j (z) = (z - p_j) / 2 + p_j$, where $p_j$ are the vertices of an equilateral triangle. Writing $\alb = \set{0,1,2}$, it has the natural \pcf\ \s-s structure $\sss_0 = (K, \alb, \set{f_j}_{j = 0}^2)$. The associated surjection $\pi: \alb^{-\omega} \to K$ induces the equivalence relations
\[
\recur{0} 1 w \sim \recur{1} 0 w, \quad \recur{1} 2 w \sim \recur{2} 1 w, \quad \text{and} \quad \recur{2} 0 w \sim \recur{0} 2 w
\]
for all $w\in \alb^*$. However, since $\recur{2}0 \sim \recur{0}2$ but $\recur{2} \not\sim \recur{0}$, the shift map $\shift$ is not well defined, and so $\sss_0$ is not a \s-s structure on the limit space of a contracting group.

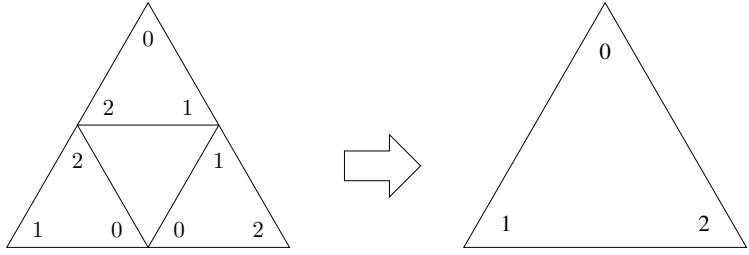
\begin{figure}
\resizebox{.8\columnwidth}{!}{\begin{tikzpicture}
\foreach \a in {0,...,2}{
\node[regular polygon, regular polygon sides=3, minimum size=2.5cm, draw] at ($(\a*120+90:1.25cm)$) {};
\node[regular polygon, regular polygon sides=5, minimum size=1.5cm] at ($(0*120+90:1.25cm)+({-1}*\a*120+90:.7cm)$) {$\a$};
\node[regular polygon, regular polygon sides=5, minimum size=1.5cm] at ($(1*120+90:1.25cm)+({-1}*\a*120+90+2*120:.7cm)$) {$\a$};
\node[regular polygon, regular polygon sides=5, minimum size=1.5cm] at ($(2*120+90:1.25cm)+({-1}*\a*120+90+4*120:.7cm)$) {$\a$};
}

\node[single arrow, draw, right] at (0:3cm) {\phantom{woot}};

\node[regular polygon, regular polygon sides=3, minimum size=5cm, draw] at (0:7cm) {};

\foreach \c in {0,...,2}{
\node[regular polygon, regular polygon sides=3, minimum size=1.5cm] at ($(0:7cm)+(\c*120+90:1.75cm)$) {\c};
}
\end{tikzpicture}}
\caption{Self-covering of Sierpi\'nski gakset in Example~\ref{eg:sierp_gasket}}
\label{fig:sierp_gasket_covering}
\end{figure}

On the other hand, there is an alternate \pcf\ \s-s structure on the Sierpi\'nski gasket, given by $\sss = (K, \alb, \set{F_j}_{j = 0}^2)$, where $F_j = r_j \circ f_j$, and $r_j$ is the reflection about the axis of symmetry through $p_j$. This can be described by the self-covering depicted in Figure~\ref{fig:sierp_gasket_covering}. With this \s-s structure, the induced equivalence relations are
\[
\recur{0} 1 w \sim \recur{0} 2 w, \quad \recur{1} 2 w \sim \recur{1} 0 w, \quad \text{and} \quad \recur{2} 0 w \sim \recur{2} 1 w \text{,}
\]
which are the relations describing the asymptotic equivalence of the 3-peg Hanoi Towers Group $G$ \cite{GS06, GS08}. In particular, $\sss$ is the natural \s-s structure on $\lims$, where $G$ is the group generated by the elements
\begin{gather*}
a_{01} = (01) (1, 1, a_{01}) \text{,}\\
a_{12} = (12) (a_{12}, 1, 1) \text{,}\\
a_{20} = (02) (1, a_{02}, 1) \text{,}
\end{gather*}
acting on the rooted tree $\alb^*$. The Moore diagram of the nucleus of the group is given in Figure~\ref{fig:hanoi_towers}. It is interesting to note that this is exactly the group $G_\sss$ that would result from our construction described in Section~\ref{sec:construction}.
\end{eg}

\begin{eg}
\label{eg:nonstrict}
This example highlights the fact that the \pcf\ condition and Condition~\ref{cond:quasi-monocarpic} together still do not imply the strictly \pcf\ condition, by describing a group that satisfies the former conditions but not the latter.

We consider the group $G$ whose nucleus is illustrated in Figure~\ref{fig:nonstrict}. This group is contracting, \pcf, and satisfies Condition~\ref{cond:quasi-monocarpic}, but is not strictly \pcfp

\begin{figure}
\resizebox{.8\columnwidth}{!}{\begin{tikzpicture}[->, shorten >=1pt, node distance=1.5cm]
\node[state] at (0,0)           (1)  {$1$};
\node[state] [left=of 1]        (h)  {$h$};
\node[state] [right=of 1]       (a)  {$a$};
\node[state] [above right=of a] (g)  {$g$};
\node[state] [below right=of a] (gh) {$gh$};
\path
(h)  edge [bend left=25]        node [above] {$(0,1)$} (1)
     edge [bend right=25]       node [below] {$(1,0)$} (1)
     edge [loop left]           node [left]  {$(2,2)$} ()
(a)  edge [bend left=40]        node [above] {$(0,1)$} (1)
     edge [bend right=40]       node [above] {$(1,0)$} (1)
     edge [bend left=0]         node [above] {$(2,2)$} (1)
(g)  edge [bend left=25]        node [right] {$(0,1)$} (a)
     edge [bend right=25]       node [left]  {$(1,0)$} (a)
     edge [in=30,out=60,loop]   node [right] {$(2,2)$} ()
(gh) edge [bend left=25]        node [left]  {$(0,0)$} (a)
     edge [bend right=25]       node [right] {$(1,1)$} (a)
     edge [in=330,out=300,loop] node [right] {$(2,2)$} ();
\end{tikzpicture}}
\caption{Nucleus of the group $G$ in Example~\ref{eg:nonstrict}}
\label{fig:nonstrict}
\end{figure}
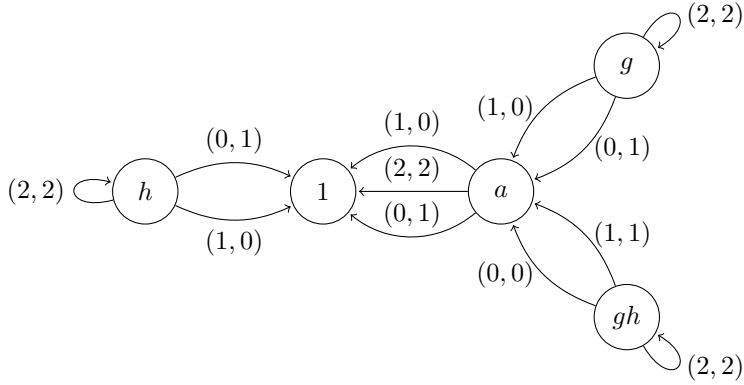

$G$ is generated by
\begin{gather*}
a = (01) (1, 1, 1) \text{,}\\
g = (01) (a, a, g) \text{,}\\
h = (01) (1, 1, h) \text{,}
\end{gather*}
and acts on $\alb = \set{0,1,2}$. Since all of its generators are bounded, $G$ is bounded and thus contracting and \pcf, with the nucleus being $\nucl = \set{1, a, g, h, gh}$. From the diagram, we see that $G$ satisfies Condition~\ref{cond:quasi-monocarpic}, and therefore there exists a \s-s structure on its limit space $\lims$. However, since $g$ can change more than one letter in a give word, this group is not strictly \pcfp
\end{eg}

\begin{eg}
\label{eg:alt_grigorchuk}
This is a straightforward example that illustrates our construction of a contracting group with a given \s-s structure.

Consider the unit interval $I = [0,1]$. Set $F_0 (x) = - (1/2) x + 1/2$, and $F_1 (x) = (1/2) x + 1/2$. Then $\sss = (I, \set{0, 1}, \set{F_i}_{i \in \set{0, 1}})$ is a \pcf\ \s-s structure on $I$. Notice that this \s-s structure is naturally isomorphic to a \s-s structure on the Koch curve. The critical set is given by $\crits = \set{\recur{1}00, \recur{1}01}$, and its image in $I$ is $\pi (\crits) = \set{1/2}$.

There exists a well-defined continuous induced shift map $\shift: I \to I$, defined by $\shift \circ \pi = \pi \circ \sigma$; we can write it explicitly as
\[
\shift (x) =
\begin{cases}
1-2x \quad \text{if } x \in [0, \frac{1}{2}] \text{,}\\
2x-1 \quad \text{if } x \in [\frac{1}{2}, 1] \text{.}
\end{cases}
\]
By Theorem~\ref{thm:construct_iff_s_exists}, the existence of $\shift$ guarantees the success of the construction of a contracting group $G_\sss$. There is only one entry in the list of equivalence classes, namely
\[
\pi^{-1} \left( \frac{1}{2} \right) = \recur{1} S_2 S_1 \text{,}
\]
where $S_2 = \set{0}$ and $S_1 = \set{0,1}$. We define
\begin{align*}
g_1 & = (01) (1, 1) \text{,}\\
g_2 & = \hspace{6.2mm} (g_1, g_2) \text{,}
\end{align*}
where we have suppressed the subscript $1/2$. Define $G_\sss = \gen{g_1, g_2}$. Notice that
\begin{align*}
g_1^2 & = (1, 1) = 1 \text{,}\\
g_2^2 & = (g_1^2, g_2^2) = (1, g_2^2) = 1 \text{,}
\end{align*}
so both generators are of order $2$. Moreover,
\begin{align*}
g_1 g_2 & = (01) (g_1, g_2) \text{,}\\
g_2 g_1 & = (01) (g_2, g_1) \text{,}
\end{align*}
so $G_\sss$ is the infinite dihedral group, and we see that the nucleus $\nucl = \set{1, g_1, g_2}$. Figure~\ref{fig:alt_grigorchuk} shows the Moore diagram of the nucleus. It can easily be seen that $G_\sss$ is in fact strictly \pcfp

\begin{figure}
\begin{tikzpicture}[->, shorten >=1pt, node distance=1.4cm]
\node[state]               (g2) {$g_2$};
\node[state] [right=of g2] (g1) {$g_1$};
\node[state] [right=of g1] (1)  {$1$};
\path
(g1) edge [bend left=28]  node [above] {$(0,1)$} (1)
     edge [bend right=28] node [below] {$(1,0)$} (1)
(g2) edge                 node [below] {$(0,0)$} (g1)
     edge [loop left]     node [left]  {$(1,1)$} (g2);
\end{tikzpicture}
\caption{Nucleus of the group $G_\sss$ in Example~\ref{eg:alt_grigorchuk}.}
\label{fig:alt_grigorchuk}
\end{figure}

The limit space of the Grigorchuk group (introduced and discussed in \cite{Gri80, Gri84}) has the same \s-s structure as $\sss$. The Grigorchuk group is also strictly \pcfp The nucleus of the Grigorchuk group is different from the nucleus of $G_\sss$, and so they are not isomorphic. In other words, it is possible for two \pcf\ groups satisfying Condition~\ref{cond:quasi-monocarpic} that are not isomorphic to each other to have limit spaces with isomorphic \s-s structures.

It has been shown that there are countably many groups whose limit space admits a \s-s structure isomorphic to $\sss$; for a classification of all such groups, see \cite{Nek03, Sun07}.
\end{eg}

\begin{eg}[Pentakun]
\label{eg:pentakun}
The pentakun, as described in \cite{Kig01}, is the unique non-empty compact space $K \subset \B{C}$ that is invariant under the injections
\[
f_k(z) = \frac{3 - \sqrt{5}}{2} (z - p_k) + p_k \text{,} \quad \text{where} \quad p_k = e^{2 \pi i k / 5}. 
\]
Figure~\ref{fig:sierp_gasket_pentakun} gives the (rotated) picture of the pentakun.

Identifying $\alb = \set{0, 1, 2, 3, 4}$, the natural \pcf\ \s-s structure is given by $\sss_0 = (K, \alb, \set{f_j}_{j = 0}^4)$. The equivalence classes induced by $\pi$ are
\[
\recur{2} 0 w \sim \recur{4} 1 w, \quad \recur{3} 1 w \sim \recur{0} 2 w,\quad \recur{4} 2 w \sim \recur{1} 3 w, \quad \recur{0} 3 w \sim \recur{2} 4 w, \quad \text{and} \quad \recur{1} 4 w \sim \recur{3} 0 w
\]
for all $w \in \alb^*$.

As with the Sierpi\'nski gasket, the shift map is not defined for this \s-s structure, and $\sss_0$ is not the \s-s structure on a limit space. However, like the Sierpi\'nski gasket, there is a modified \pcf\ \s-s structure that can be achieved as the \s-s structure on the limit space of a contracting group.

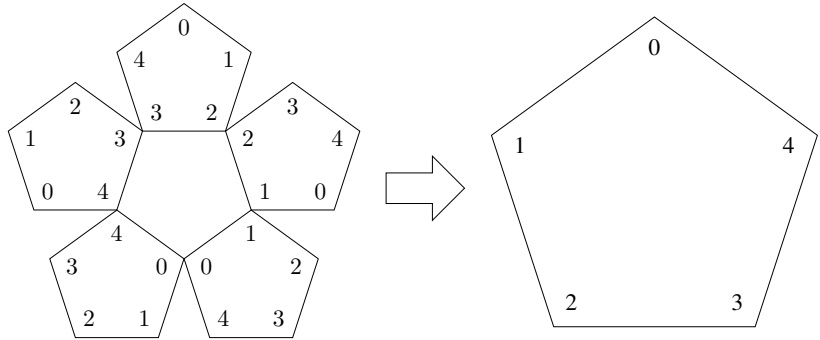
\begin{figure}
\resizebox{.9\columnwidth}{!}{\begin{tikzpicture}
\foreach \a in {0,...,4}{
\node[regular polygon, regular polygon sides=5, minimum size=2.1cm, draw] at ($(\a*72+90:1.7cm)$) {};
\node[regular polygon, regular polygon sides=5, minimum size=1.5cm] at ($(0*72+90:1.7cm)+({-1}\a*72+90:.7cm)$) {$\a$};
\node[regular polygon, regular polygon sides=5, minimum size=1.5cm] at ($(1*72+90:1.7cm)+({-1}*\a*72+90+2*72:.7cm)$) {$\a$};
\node[regular polygon, regular polygon sides=5, minimum size=1.5cm] at ($(2*72+90:1.7cm)+({-1}*\a*72+90+4*72:.7cm)$) {$\a$};
\node[regular polygon, regular polygon sides=5, minimum size=1.5cm] at ($(3*72+90:1.7cm)+({-1}*\a*72+90+6*72:.7cm)$) {$\a$};
\node[regular polygon, regular polygon sides=5, minimum size=1.5cm] at ($(4*72+90:1.7cm)+({-1}*\a*72+90+8*72:.7cm)$) {$\a$};
}

\node [single arrow, draw, right] at (0:3cm) {\phantom{woot}};

\node[regular polygon, regular polygon sides=5, minimum size=5.1cm,draw] at (0:7cm) {};

\foreach \c in {0,...,4}{
\node[regular polygon, regular polygon sides=5, minimum size=1.4cm] at ($(0:7cm)+(\c*72+90:2.1cm)$) {\c};
}
\end{tikzpicture}}
\caption{Self-covering of pentakun in Example~\ref{eg:pentakun}}
\label{fig:pentakun_covering}
\end{figure}

Consider $\sss = (K, \alb, \set{F_j}_{j=0}^4)$, where $F_j = r_j \circ f_j$, and $r_j$ is the reflection about the line joining $p_j$ with the origin, i.e.\ the axis of symmetry through $p_j$ of the pentagon formed by $\set{p_j}_{j=0}^4$. The corresponding self-covering is depicted in Figure~\ref{fig:pentakun_covering}. Here the equivalence classes are of the form $\recur{k} S_k w$ where $k \in \alb$ and $S_k = \set{k-2 \mod 5, \; k+2 \mod 5}$.

\begin{figure}
\resizebox{.6\columnwidth}{!}{\begin{tikzpicture}[->, shorten >=1pt, node distance=2cm]
\node[state] at (0,0)     (1)  {$1$};
\node[state] at (90:2.5)  (a0) {$a_0$};
\node[state] at (162:2.5) (a1) {$a_1$};
\node[state] at (234:2.5) (a2) {$a_2$};
\node[state] at (306:2.5) (a3) {$a_3$};
\node[state] at (378:2.5) (a4) {$a_4$};
\path
(a0) edge [loop above]          node [above]                  {$(0,0)$} ()
     edge [bend right=20]       node [near start, left]       {$(2,3)$} (1)
     edge [bend left=20]        node [near start, right]      {$(3,2)$} (1)
(a1) edge [in=177,out=147,loop] node [above]                  {$(1,1)$} ()
     edge [bend right=20]       node [very near start, below] {$(3,4)$} (1)
     edge [bend left=20]        node [near start, above]      {$(4,3)$} (1)
(a2) edge [in=249,out=219,loop] node [below]                  {$(2,2)$} ()
     edge [bend right=20]       node [near start, below]      {$(4,0)$} (1)
     edge [bend left=20]        node [left]                   {$(0,4)$} (1)
(a3) edge [in=321,out=291,loop] node [below]                  {$(3,3)$} ()
     edge [bend right=20]       node [right]                  {$(0,1)$} (1)
     edge [bend left=20]        node [near start, below]      {$(1,0)$} (1)
(a4) edge [in=370,out=400,loop] node [above]                  {$(4,4)$} ()
     edge [bend right=20]       node [near start, above]      {$(1,2)$} (1)
     edge [bend left=20]        node [very near start, below] {$(2,1)$} (1);
\end{tikzpicture}}
\caption{Nucleus of the group $G_\sss$ in Example~\ref{eg:pentakun}}
\label{fig:pentakun_group}
\end{figure}

Our construction from Section~\ref{sec:construction} yields the \pcf\ group $G_\sss$ generated by
\begin{gather*}
a_0 = (23) (a_0, 1, 1, 1, 1) \text{,}\\
a_1 = (34) (1, a_1, 1, 1, 1) \text{,}\\
a_2 = (40) (1, 1, a_2, 1, 1) \text{,}\\
a_3 = (01) (1, 1, 1, a_3, 1) \text{,}\\
a_4 = (12) (1, 1, 1, 1, a_4) \text{,}
\end{gather*}
so that $\sss$ is a \s-s structure on the limit space of $G_\sss$. The Moore diagram of the nucleus of $G_\sss$ is shown in Figure~\ref{fig:pentakun_group}.

It is easy to perform an analogous construction for all $n$-kuns where $n$ is odd.
\end{eg}

\begin{eg}[Hexakun and Linstr\o m Snowflake]
\label{eg:hexakun_snowflake}
In Example~\ref{eg:pentakun} we showed how a self-covering could be constructed on the pentakun that could be taken to be the shift map required for the construction in Section~\ref{sec:construction}. The current example shows the way to construct a self-covering for the hexakun, a fractal analogous to the pentakun but constructed instead from a hexagon. We shall also discuss why no self-covering can be constructed for the Linstr{\o}m snowflake, a nested fractal which is a variation of the hexakun.

Similar to the pentakun, the hexakun is typically constructed (e.g.\ in \cite{Kig01}) as the unique non-empty compact space $K \subset \B{C}$ invariant under the injections
\[
f_k(z) = \frac{1}{3}(z- p_k) + p_k \text{,} \quad \text{where} \quad p_k = e^{\pi i k/3} \text{.}
\]
Writing $\alb = \set{0,1,2,3,4,5}$, we see that $\sss_0 = (K, \alb, \set{f_j}_{j=0}^5)$ is the usual \s-s structure. As with the Sierpi\'nski gasket and the pentakun, this \s-s structure does not admit a shift map, and so we have to choose another \s-s structure.

\begin{figure}
\resizebox{.9\columnwidth}{!}{\begin{tikzpicture}
\foreach \a in {0,...,5}{
\node[regular polygon, regular polygon sides=6, minimum size=1.5cm, draw] at ($(\a*60:1.5cm)$) {};
\node[regular polygon, regular polygon sides=6, minimum size=1cm] at
($(0*60:1.5cm)+(\a*60:.5cm)$) {\small$\a$};
\node[regular polygon, regular polygon sides=6, minimum size=1cm] at
($(1*60:1.5cm)+(60+{-1}*\a*60:.5cm)$) {\small$\a$};
\node[regular polygon, regular polygon sides=6, minimum size=1cm] at
($(2*60:1.5cm)+(60*2+\a*60:.5cm)$) {\small$\a$};
\node[regular polygon, regular polygon sides=6, minimum size=1cm] at
($(3*60:1.5cm)+(60*3+{-1}*\a*60:.5cm)$) {\small$\a$};
\node[regular polygon, regular polygon sides=6, minimum size=1cm] at
($(4*60:1.5cm)+(60*4+\a*60:.5cm)$) {\small$\a$};
\node[regular polygon, regular polygon sides=6, minimum size=1cm] at
($(5*60:1.5cm)+(60*5+{-1}*\a*60:.5cm)$) {\small$\a$};
}

\node [single arrow, draw, right] at (0:2.7cm) {\phantom{woot}};

\node [regular polygon, regular polygon sides=6, minimum size=4.5cm, draw] at (0:6.5cm) {};

\foreach \c in {0,...,5}{
\node[regular polygon, regular polygon sides=6, minimum size=1cm] at
($(0:6.5cm)+(\c*60:1.9cm)$) {\small\c};
}
\end{tikzpicture}}
\caption{Self-covering of hexakun in Example~\ref{eg:hexakun_snowflake}}
\label{fig:hexakun_covering}
\end{figure}
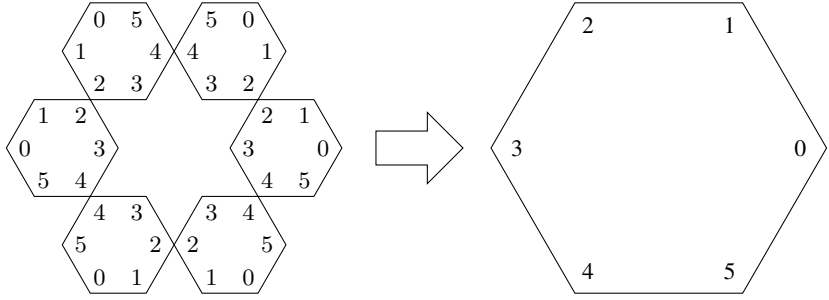

This fractal is a set of six copies of itself, each with a ``corner'' at the point $p_k$, and joined to two adjacent copies at the corner, two corners away from $p_k$. Our self-covering can be thought of as folding the fractal in half along the $y$-axis, so the cells in the left half-plain land on their reflections in the right. For the 3 cells on the right, we fold the upper and lower cells onto the cell containing $p_0$. Finally, we rescale the $p_0$ cell using the map $f_0^{-1}$.

Formally, if we let $\phi_1(z) = e^{\pi i/3}\recur{z}$, $\phi_2(z) = e^{2\pi i/3}z$, $\phi_3(z) = -\recur{z}$, $\phi_4(z) = e^{4\pi i/3}z$, and $\phi_5(z) = e^{5\pi i/3}\recur{z}$ (here $\recur{z}$ is the complex conjugate), then our modified \s-s structure is $\sss = (K, \alb, \set{F_j}_{j=0}^5)$ where $F_0 = f_0$ and $F_j = \phi_j \circ f_j$ for $j = 1, 2, 3, 4, 5$.

We now find the critical set $\crits$ and post-critical set $\pcrits$ of this new \s-s structure $\sss$. Noticing that the fixed point of $F_0$ is still $p_0$, we see that $\pi (\recur{0}) = p_0$. Also, $F_j (p_0) = p_j$, and so $\pi (\recur{0} j) = p_j$ for $j = 1, 2, 3, 4, 5$. Since cells are only joined at the corners, this is enough to give us the addresses for the entire critical set; thus, $\crits = \set{\recur{0} 2 j, \recur{0} 4 j \; | \; j \in \alb}$. Notice that in this \s-s structure, only the points $p_2$ and $p_4$ are mapped to boundary points. The post-critical set is then given by $\pcrits = \set{\recur{0}, \recur{0} 2,\recur{0} 4}$.

More precisely, we examine the \s-s structure and write down the equivalence classes as follows:
\begin{gather*}
\recur{0} 2 0 w \sim \recur{0} 2 1 w, \quad \recur{0} 2 2 w \sim \recur{0} 2 3 w, \quad \recur{0} 2 4 w \sim \recur{0} 2 5 w,\\
\recur{0} 4 0 w \sim \recur{0} 4 5 w, \quad \recur{0} 4 1 w \sim \recur{0} 4 2 w, \quad \recur{0} 4 3 w \sim \recur{0} 4 4 w.
\end{gather*}
Applying our construction from Section~\ref{sec:construction}, we get a group $G_\sss$ generated by the nine elements with wreath recursions
\[
a_{01} = (01), \quad a_{23} = (23), \quad a_{45} = (45),
\]
and
\begin{gather*}
b_0 = (b_0, 1, a_{01}, 1, 1, 1), \quad b_1 = (b_1, 1, a_{23}, 1, 1, 1), \quad b_2 = (b_2, 1, a_{45}, 1, 1, 1),\\
b_3 = (b_3, 1, 1, 1, a_{01}, 1), \quad b_4 = (b_4, 1, 1, 1, a_{23}, 1), \quad b_5 = (b_5, 1, 1, 1, a_{45}, 1).
\end{gather*}
Notice that according to our construction, there are six elements in the image $\pi (\crits)$ of the critical set, and each of these corresponds to 2 generators, and so one may expect to obtain twelve generators. However, upon closer examination, we see that, for example, $\pi (\recur{0} 2 0) = \pi(\recur{0} 2 1)$ and $\pi (\recur{0} 4 0) = \pi (\recur{0} 4 1)$ both give rise to the generator with wreath recursion $(01)$; thus, we see that three generators are redundant, and so $G_\sss$ is generated by nine elements.

\begin{figure}
\begin{center}
\includegraphics[width=0.47\linewidth]{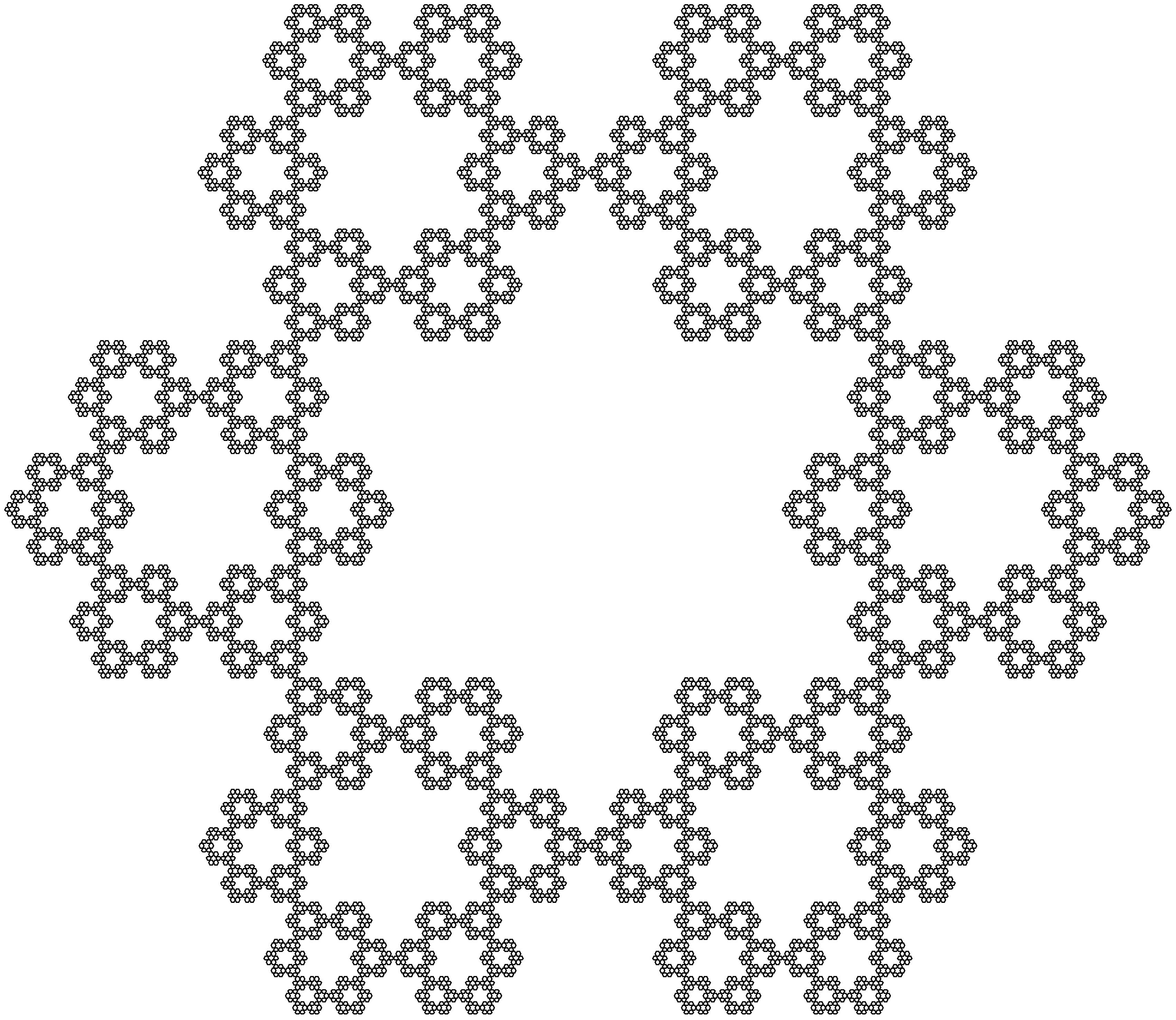}\hspace{12pt}
\includegraphics[width=0.47\linewidth]{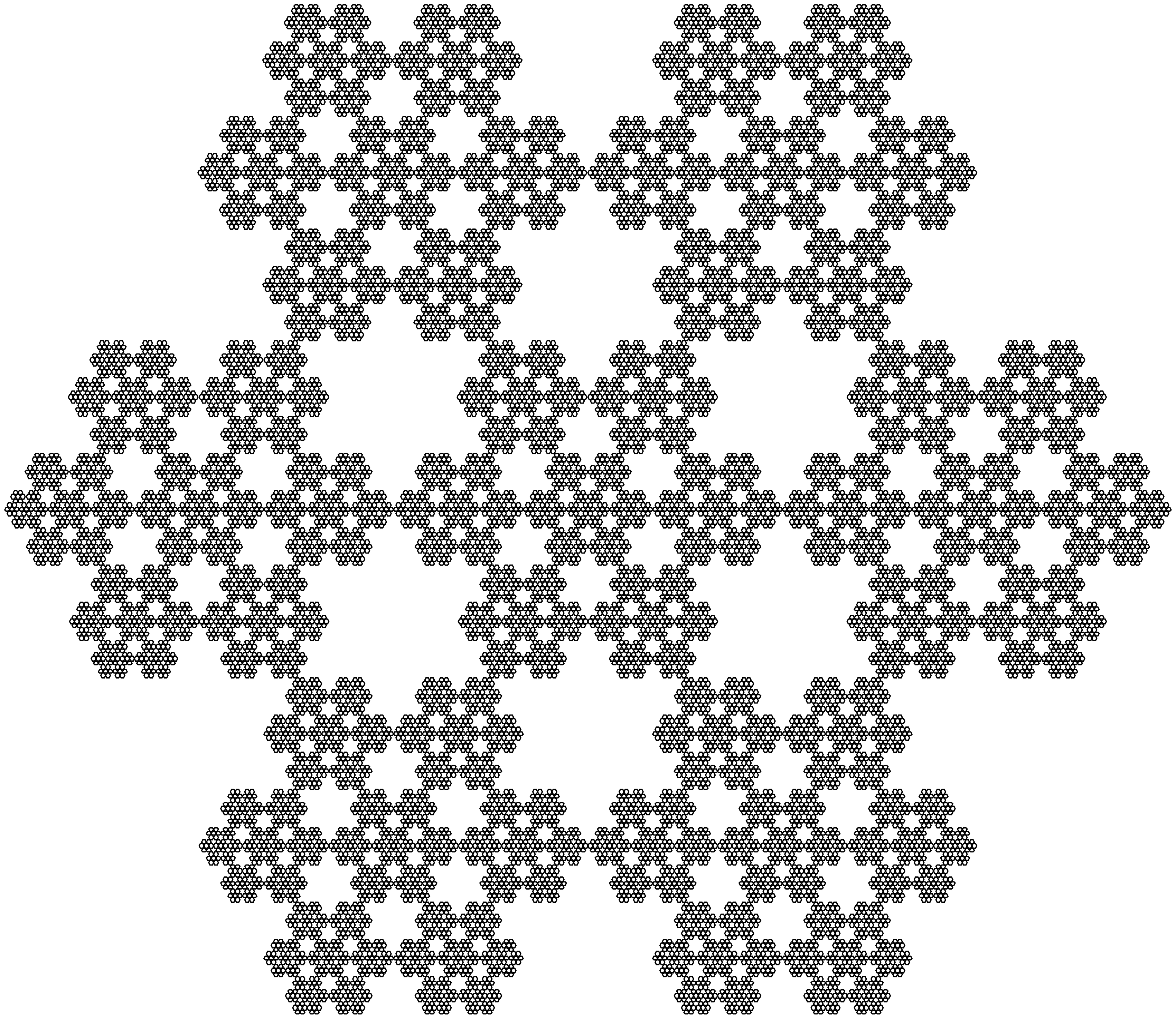}
\end{center}
\caption{Hexakun and Lindstr{\o}m snowflake}
\label{fig:hexakun_snowflake}
\end{figure}

We now turn to the Linstr{\o}m snowflake, which is a variation on the hexakun. It is the unique non-empty compact space $L \subset \B{C}$ invariant under the injections $f_0, \ldots, f_6$, where $f_j$ are the same as above for $0 \leq j \leq 5$, and $f_6 (z) = z / 3$. This fractal is like the hexakun, but with a scaled copy of itself inserted into the center. Cells of the snowflake still only intersect at the ``corners;'' in particular, $f_j(L) \cap f_k(L) = \set{f_j (p_n) \; | \; f_j (p_n) = f_k (p_m) \text{ for some } m}$ contains at most one element.

Suppose now that there exists some \s-s structure $\sss' = (L, \alb, \set{g_j}_{j=0}^6)$ on $L$ that has a shift map $\sigma$. Without loss of generality, we can assume that $g_6 (L)$ is the first-level scaled copy of $L$ in the center. Notice that $g_6 (L)$ intersects every other cell at one point, so $g_6 (p_0) = g_j (p_0)$ for some $j$, and $g_6 (p_1) = g_k (p_1)$ for some other $k$. Then $g_j (L)$ must intersect $g_k (L)$ at one point $p$, such that $g_j^{-1} (p) = g_k^{-1} (p) = \sigma (p)$ is a boundary point adjacent to both $p_0$ and $p_1$. Since no such boundary point exists, we have arrived at a contradiction. Therefore, there exists no \s-s structure on $L$ that admits a shift map, and so the Linstr{\o}m snowflake cannot arise as the limit space of a contracting group.
\end{eg}

\bibliography{Fractals}{}
\bibliographystyle{amsalpha}
%
%
\end{document}